\documentclass[preprint,12pt]{elsarticle}
\usepackage[margin=1in]{geometry}
\usepackage[utf8]{inputenc}
\usepackage{graphicx}% Include figure files
\usepackage{multirow}
\usepackage{dcolumn}% Align table columns on decimal point
\usepackage{bm}
\usepackage{color}      % color content

\usepackage{amsmath}
\usepackage{amssymb}
\begin{document}
\begin{frontmatter}

\title{Fourier neural operator based fluid-structure interaction for predicting the vesicle dynamics}
\author[1,2,3]{Wang Xiao}
\author[1,2,3]{Ting Gao \corref{cor1}}
\author[4]{Kai Liu}
\author[3,5,6]{Jinqiao Duan}
\author[1,2,3]{Meng Zhao \corref{cor2}}
\address[1]{School of Mathematics and Statistics, Huazhong University of Science and Technology, Wuhan 430074, China}
\address[2]{Center for Mathematical Sciences, Huazhong University of Science and Technology, Wuhan 430074, China}
\address[3]{Steklov-Wuhan Institute for Mathematical Exploration, Huazhong University of Science and Technology, Wuhan 430074, China}
\address[4]{College of Education for the Future, Beijing Normal University, Zhuhai 519087, China}
\address[5]{Department of Mathematics and Department of Physics, Great Bay University, Dongguan, Guangdong 523000, China}
\address[6]{Dongguan Key Laboratory for Data Science and Intelligent Medicine, Dongguan, Guangdong 523000, China}

\cortext[cor1]{Corresponding author: tgao0716@hust.edu.cn}
\cortext[cor2]{Corresponding author: mzhao9@hust.edu.cn}

\begin{abstract}
% Fluid-structure interaction analysis has high computing cost when using computational fluid dynamics. To overcome this challenge, Fourier neural operator (FNO) based fluid-structure interaction (FSI) solver is proposed. This algorithm is then applied to the interaction between an elongation flow and a vesicle. A Fourier neural operator is trained to predict incompressible flow and the dynamics of a vesicle is modeled by an integral equation using the immersed boundary method. The results show that the FNO-based FSI solver is capable of capturing the variations of the fluid and the vesicle. Moreover, ‘sequence to one’ prediction by FNO achieves higher accuracy in short-term predictions, while ‘sequence to sequence’ of FNO has higher accuracy in long-term predictions. We also provide detailed validation results on the performance comparison between interpolation and extrapolation on various scenarios.

Solving complex fluid-structure interaction (FSI) problems, characterized by nonlinear partial differential equations, is crucial in various scientific and engineering applications. Traditional computational fluid dynamics (CFD) solvers are insufficient to meet the growing requirements for large-scale and long-period simulations. Fortunately, the rapid advancement in neural networks, especially neural operator learning mappings between function spaces, has introduced novel approaches to tackle these challenges via data-driven modeling. In this paper, we propose a Fourier neural operator-based fluid-structure interaction solver (FNO-based FSI solver) for efficient simulation of FSI problems, where the solid solver based on the finite difference method is seamlessly integrated with the Fourier neural operator to predict incompressible flow using the immersed boundary method. We analyze the performance of the FNO-based FSI solver in the following three situations: training data with or without the steady state, training method with one-step label or multi-step labels, and prediction in interpolation or extrapolation. We find that the best performance for interpolation is achieved by training the operator with multi-step labels using steady-state data. Finally, we train the FNO-based FSI solver using this optimal training method and apply it to vesicle dynamics. The results show that the FNO-based FSI solver is capable of capturing the variations in the fluid and the vesicle.

\end{abstract}


\end{frontmatter}

\section{Introduction}
The fluid-structure interaction (FSI) problem is a topic of great interest in both academia and engineering \cite{2021Fluid, 2017The}. A detailed analysis of the FSI system is essential during the initial stage of the system design. The advancement of computer-aided engineering (CAE) has led to significant progress in numerical methods for FSI problems. Compared to experimental approaches, numerical simulations significantly reduce costs. However, the computational fluid dynamics (CFD) simulations are hindered by the large number of mesh grids, resulting in slow simulation speeds. The simulation of FSI problems often involves solving the structure solver after waiting for the CFD solver to provide data. Therefore, improving the simulation speed of CFD is crucial for accelerating FSI simulations.

Vesicles are sacs with fluid inside and have a nearly inextensible lipid-bilayer boundary. The dynamics of vesicles in an elongation flow is a well-studied fluid-structure interaction phenomenon. These highly nonlinear relaxation dynamics have received significant attention in numerous theoretical \cite{seifert_fluid_1999, WOS:000222502700094, WOS:000234758100106, WOS:000243893700049, abreu_noisy_2013, 2014Membrane, turitsyn_wrinkling_2008}, experimental \cite{WOS:000071043500110, WOS:000173849300065, WOS:000234118500080, WOS:000234969300067, WOS:000237124400002, deschamps_dynamics_2009, kumar_double-mode_2020, 2014Complex, kantsler_vesicle_2007}, and numerical \cite{turitsyn_wrinkling_2008, knorr_wrinkling_2010, levant_amplification_2012, liu_nonlinear_2014, liu_wrinkling_2016, WOS:A1996VM89000050, WOS:000188946500050, WOS:000225785200075, narsimhan_pearling_2015, yazdani_three-dimensional_2012, liu2019hydrodynamics, liu2019efficient, PhysRevE.107.035103} studies. In 2013, Liu {\it et al.} developed a spectral accurate boundary integral method to simulate the nonlinear evolution of the vesicle in two dimensions \cite{liu_nonlinear_2014}. Recently, Xiao {\it et al.} \cite{PhysRevE.107.035103} numerically investigated the interaction between the three-dimensional (3D) vesicle and the elongation flow using the immersed boundary method (IBM). 

% However, despite their use of GPU acceleration, the computational cost of the IBM is high. They found that the IBM solving fluid consumes up to 90\% of computing time. 

% In order to decrease the computational demands of FSI simulations such as dynamics of vesicle, we propose to use a neural operator to predict fluid flow instead of the IBM solver. More specifically, we train a Fourier neural operator (FNO) to predict fluid flow with moving boundary and couple it with a classical structure solver.

The endeavor of replacing traditional solvers with a network for predicting the fluid field and coupling it with a structural dynamic solver to solve FSI problems is becoming increasingly popular. Halder {\it et al.} \cite{halder2020deep} used recurrent neural networks (RNNs) to predict forces at the next moment based on aileron motion at this moment of an airfoil–gust and aero elastic interaction system. Han {\it et al.} \cite{han2022deep} constructed a DNN-based reduced-order model for FSI to enable the rapid prediction of FSI systems. Bublík {\it et al.} \cite{bublik2023neural} introduced an FSI solver that uses a CNN to predict fluid flow. Mazhar {\it et al.} \cite{mazhar2023novel} presented a novel interface technique employing Artificial Neural Networks (ANNs) for efficient data transfer between fluid and solid domains. 

To predict fluid flow at high speed, neural operators have been rapidly developing. E {\it et al.} \cite{yu2018deep} proposed a Deep Ritz Method for numerically solving variational problems, particularly those arising from partial differential equations. DeepGreen: deep learning of Green's functions is proposed for nonlinear boundary value problems \cite{gin2021deepgreen}. A recent approach \cite{li2020neural, li2020fourier} introduced a Fourier neural operator based on the Fourier transform that solves the problem using functional parametric dependence and learns from infinite-dimensional mappings directly. DeepONets \cite{lu2021learning} provide a solution to the numerical approximations using deep learning. \cite{CHEN2022110996} proposed Meta multigrid networks (Meta-MgNet) for solving parameterized partial differential equations. A multi-pole graph neural network (MGNO) was proposed in \cite{li2020multipole}. Similarly, U-FNO \cite{wen2022u} enhanced the FNO to work for multi-phase flow by adding a U-Net within the FNO. In this context, the FNO has made a significant breakthrough and demonstrated state-of-the-art performance. Recent scholars have conducted theoretical research on the FNO. \cite{kovachki2021universal} proved that FNOs are universal, in the sense that they can approximate any continuous operator with the desired accuracy. \cite{lanthaler2023nonlocal} introduced nonlocal neural operators (NNOs) which unifies the analysis of a wide range of neural operator architectures. To decrease the computational demands of FSI simulations, such as the dynamics of a vesicle, we use a neural operator to predict fluid flow instead of the IBM solver. More specifically, we train a Fourier neural operator (FNO) to predict fluid flow with a moving boundary and couple it with a classical structure solver.

\textbf{Our Contributions.} We introduce a FNO-based FSI solver to enable the rapid prediction of FSI systems. 

\begin{itemize}
    \item We couple the FNO and the finite difference method for solving the Navier-Stokes equation of the fluid and the integral equation of the vesicle.
    % We divide the problem into two subdomains - fluid domain and solid domain. For the fluid domain, we replace the classical method with a FNO. For the solid domain, we model the dynamics of the structure by an integral equation equation using IBM. This method exploits the Dirac-delta function to transfer the velocity from fluid to the structure. We get the position of the structure by solving the integral equation using a finite different method. 
    \item We apply the FNO-based FSI solver to vesicle dynamics. This solver produces accurate results with high computational speed. It is an order of magnitude more efficient than IBM.
    % We test the FNO-based FSI solver on a problem with vesicle dynamic. In particular, we solve the interaction between an elongation flow and a vesicle sacs with fluid inside and have a nearly inextensible lipid-bilayer boundary. The FNO model succeed in predicting the unsteady flow field around a moving object. The predictions from 'sequence to sequence' types of FNO are used in simulating the dynamic process of the vesicle.
    \item We analyze the performance of the algorithm in the following three situations: training data with or without the steady state, training method with one-step label or multi-step labels, and prediction in interpolation or extrapolation. Using data with steady state to train operator with multi-step labels for interpolation performs best.
    % Analysis on performance of algorithm. We set four experiments and investigate the performance of different training types of FNO. We find that 'sequence to one' types of FNO achieve higher accuracy in short-term predictions, while 'sequence to sequence' types of FNO have higher accuracy in long-term predictions. For the 'sequence to one' types of FNO, one without steady-state data performs batter than another one with steady-state data in doing extrapolation.For the 'sequence to sequence' types of FNO, one with steady-state data performs batter than another one without steady-state data.   
\end{itemize}

The paper is organized as follows. In section 2, we provide the problem description of vesicle dynamics. In section 3, we introduce the FNO-based FSI solver. In section 4, we present numerical results and discussions. In section 5, we offer concluding remarks and discuss potential future work.

\section{Problem description}
Xiao {\it et al.} \cite{PhysRevE.107.035103} studied the three-dimensional dynamic of a vesicle in an elongation flow by the immersed boundary method. Here we review this problem. An external extensional flow field ${\bf u}^{\infty}({\bf x},t)$ is applied for ${\bf x} \in \Omega$. The external extensional flow $u_x^{\infty}=\dot{\gamma}x$, $u_y^{\infty}=-\dot{\gamma}y$, and $u_z^{\infty}=0$, where $\dot{\gamma}$ is the extensional rate. The total velocity field ${\bf u}({\bf x},t)$ is decomposed into the external elongational flow ${\bf u}^{\infty}$ and the induced flow ${\bf u}^{\rm ind}$ due to the presence of the vesicle \cite{liu_nonlinear_2014, PhysRevE.107.035103}, ${\bf u}({\bf x},t) = {\bf u}^{\infty}({\bf x})+{\bf u}^{\rm ind}({\bf x},t)$. The vesicle moves with the fluid with velocity ${\bf u}({\bf x},t)$. When the elongation flow is suddenly turned on, the vesicle relaxes towards a stationary state. Then, when we set the elongation flow suddenly switched from $u_x^{\infty}=\dot{\gamma}x$ and $u_y^{\infty}=-\dot{\gamma}y$ to $u_x^{\infty}=-\dot{\gamma}x$ and $u_y^{\infty}=\dot{\gamma}y$, while $u_z^{\infty}$ remains zero, the vesicle undergoes a relaxation from the stretched stationary state to another one \cite{kantsler_vesicle_2007, liu_wrinkling_2016, PhysRevE.107.035103}.

A 3D vesicle is placed within a viscous fluid and contains the same fluid internally. The fluids inside and outside the vesicle are assumed to be highly viscous and satisfy the Stokes equation \cite{liu_nonlinear_2014, PhysRevE.107.035103},
\begin{align}
&\rho \frac{\partial {\bf u}^{\rm ind}({\bf x},t)}{\partial t} = \mu \Delta {\bf u}^{\rm ind}({\bf x},t) - \nabla p + {\bf f}({\bf x},t),
\label{fluidv}\\
&\nabla \cdot {\bf u}^{\rm ind}({\bf x},t)  =  0, \label{incom} 
\end{align}
where ${\bf u}^{\rm ind}({\bf x},t)$ is the Eulerian velocity field of the fluid at the position ${\bf x} \in  {\mathbb R}^3$ and time $t \in \mathbb{R}^+$ in the fluid domain $\Omega$, $p$ is the pressure, $\mu$ is the dynamic viscosity, $\rho$ is the density of the fluid, and ${\bf f}({\bf x},t)$ is the force density acting on the fluid. 

In the typical IBM, the force density ${\bf f}({\bf x},t)$ exerted on the fluid is obtained by applying a delta function $\delta({\bf x})$ to convert the force ${\bf F}({\bf X},t)$ associated with the vesicle membrane \cite{PhysRevE.107.035103, peskin_immersed_2002, atzberger2007stochastic}. Here ${\bf X}$ is the position of the membrane. By interpolating the velocity of the nearby flow field, we establish the fluid-vesicle coupling as follows,
\begin{align}
{\bf f}({\bf x},t) &=\int_{\mathcal{S}} \delta[{\bf x}-{\bf X}(t)]{\bf F}({\bf X},t) d{\bf X}, \label{focon}  \\
{\frac{d{\bf X}(t)}{dt}} &= {\int_{\Omega}{\delta[{\bf x}-{\bf X}(t)]{\bf u}({\bf x},t) d{\bf x}}}, \label{coupling} 
\end{align}
where $\mathcal{S}$ is the membrane of vesicle. 

\section{FNO-based FSI solver}
To solve FSI problems, we divide the problem into two subdomains - the fluid domain and the solid domain. In the fluid domain, we replace the classical method with FNO, serving as a fluid dynamic solver. In the solid domain, we model the dynamics of the structure using an integral equation with IBM. This method utilizes the Dirac-delta function to transfer velocity from the fluid to the structure. We obtain the position of the structure by solving the integral equation using a finite difference method.

% The traditional FSI numerical simulation framework consists of two parts: a fluid dynamic solver and a structural dynamic solver. The FSI framework in this paper still adopts the same strategy. This FSI framework also consists of two parts: a fluid dynamic solver (DNN based fluid model) and a structural dynamic solver (structural motion control equations). The velocity of fluid is transferred between the fluid  and the structural dynamics solver.

\subsection{Review of the Fourier neural operator}
The Fourier neural operator \cite{li2020fourier, li2020neural} learns a mapping between two infinite dimensional spaces from a finite collection of observed input-output pairs. To formulate the problem, we define the domain $D \subset \mathbb{R}^d$ be a bounded and open set; $\mathcal{A}$ and $\mathcal{U}$ be separable Banach spaces of function defined on $D$ that take values in $\mathbb{R}^{d_a}$ and $\mathbb{R}^{d_u}$ respectively. $G^{\dagger}: \mathcal{A} \to \mathcal{U}$ is a non-linear map that satisfies the governing PDEs. Suppose we have observations $\{a_j, u_j\}_{j=1}^N$ where $a_j \sim u$ is an i.i.d. sequence from the probability measure supported on $\mathcal{A}$ and $u_j = G^{\dagger} (a_j)$ is possibly corrupted with noise. We aim to build an operator $G_\theta$ that learns an approximation of $G^{\dagger}$ by minimizing the following problem, 
\begin{equation}
    \mathop{\min}_{\theta} \mathbb{E} [ \left ( \frac{norm(G(a,\theta), G^{\dagger}(a))}{norm(G^{\dagger}(a))}\right)].
\end{equation}
Since $a_j \in \mathcal{A}$ and $u_j \in \mathcal{U}$ are both functions, we use $n$-point discretization $D_j = \{ x_1, ..., x_n \} \subset D$ to numerically represent $a(x)_j|_{D_j} \in \mathbb{R}^{n\times d_a}$ and $u(x)_j|_{D_j} \in \mathbb{R}^{n\times d_u}$. 

The FNO is formulated as an iterative architecture,
\begin{equation}
    a(x) \mapsto v_0 \mapsto v_1 \mapsto ... \mapsto v_E \mapsto u(x),
\end{equation}
 where $v_j$ for $j=0,1,...,E-1$ is a sequence of functions each taking values in $\mathbb{R}^{d_v}$. The input $a\in \mathcal{A}$ is initially transformed to a higher dimensional representation $v_0(x)=P(a(x))$ through the local transformation $P$, typically parameterized by a shallow fully-connected neural network. Subsequently, we perform several iterations of updates $v_t \mapsto v_{t+1}$ (defined below). The output $u(x) = Q(v_E(x))$ corresponds to the projection of $v_E$ using the local transformation $Q:\mathbb{R}^{d_v}\to \mathbb{R}^{d_u}$. In each iteration, the update $v_t \mapsto v_{t+1}$ is defined as the composition of a non-local integral operator $\mathcal{K}$ and a local, nonlinear activation function $\sigma$,
\begin{equation}
    v_{t+1}(x) := \sigma \left(Wv_t(x) + ( \mathcal{K} (a; \phi)v_t)(x)  \right), \qquad \forall x \in D \label{FNO1}
\end{equation}
where $W: \mathbb{R}^{d_v} \to \mathbb{R}^{d_v}$ is a linear transformation and $\sigma: \mathbb{R} \to \mathbb{R}$ is a non-linear activation function whose action is defined component-wise. Here $\mathcal{K}$ maps to bounded linear operators and is defined by 
\begin{equation}
    \left( \mathcal{K}(\phi)v_t \right)(x) = \mathcal{F}^{-1} \left(R \cdot (\mathcal{F}v_t)\right)(x), \qquad \forall x \in D \label{FNO3}
\end{equation}
where $\mathcal{F}$ denote the Fourier transform of a function $f: D \to \mathbb{R}^{d_v}$, $\mathcal{F}^{-1}$ is its inverse, and $R$ is the Fourier transform of a periodic function $\kappa: \bar D \to \mathbb{R}^{d_v \times d_v}$. Since we assume that $\kappa$ is periodic, we can apply a Fourier series expansion and work in the discrete modes of the Fourier transform. We first truncate the Fourier series at a maximum number of modes $k_{max}$, and then parameterize $R$ directly as a complex-valued $(k_{max}\times c\times c)$-tensor with the truncated Fourier coefficients. As a result, multiplication by the learnable weight tensor $R$ is 
\begin{equation}
    \left( R \cdot \mathcal{F}(v_t) \right)_{k,i} = \sum_{j=1}^c R_{k,i,j} \left( \mathcal{F}(v_t) \right)_{k,j}, \qquad \forall k=1,...,k_{max}, i=1,...,c.
\end{equation}
By replacing the $\mathcal{F}$ by the FFT and implementing $R$ using a direct linear parameterization, we have obtained the Fourier operator.

\subsection{The FSI solver on the problem}
For the dynamic of vesicle (\ref{fluidv}) - (\ref{coupling}), we use an Eulerian grid to discretize the flow field and a Lagrangian grid to discretize the vesicle. The flow field is a cube with an edge $L$ in all directions. The grid size $\Delta x= L/N$, where $N$ is the number of grid points along each direction. Let ${\bf{u}}_{\bf{m}}^{{\rm ind},n} = {\bf{u}}_{\bf{m}}^{{\rm ind}}(t_{n})$ be the induced velocity at the ${\bf m}^{th}$ grid point at discrete time $t_n=n\Delta t$. Here ${\bf m}=(m_1,m_2,m_3)$ is a vector with integer components. The vesicle is discretized by a triangular mesh with edge $e \approx 2\Delta x$ and a total number of vertices $M$ \cite{persson2004simple}. Let ${\bf X}_k^n$ be the position of the $k$-th mesh point at discrete time $t_n=n\Delta t$.

\begin{figure}
    \centering
    \includegraphics[width=1\linewidth]{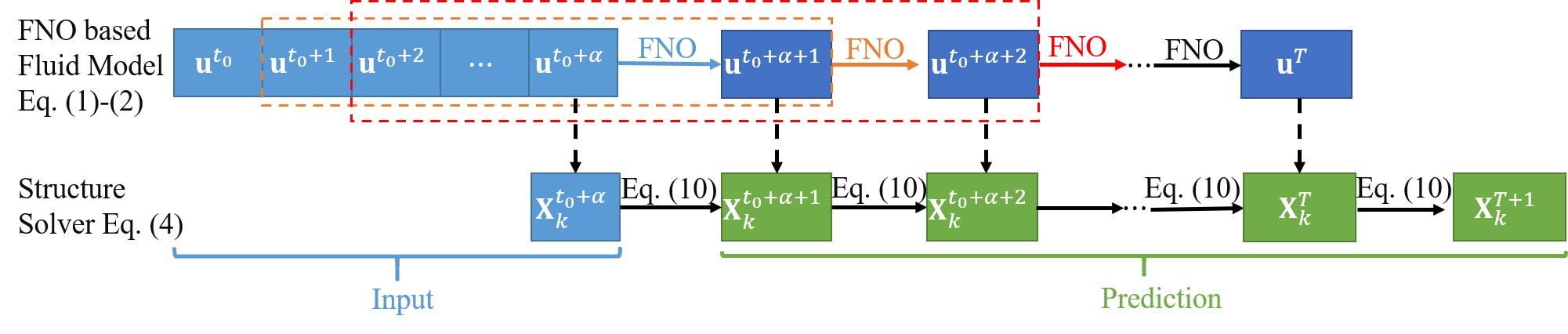}
    \caption{The architecture of the FNO-based FSI simulation framework. \textbf{Top:} The architecture of the FNO-based fluid model. The operator maps the solution at the previous $\alpha$ time steps to the next time step. For example, the data from a dashed yellow box is used to predict ${\bf u}^{t_0+\alpha+2}$, with a yellow arrow connecting the feature and label. The FNO is achieved using a combination of 3D-convolutional layers for spatial information processing and recurrent layers for temporal information propagation. \textbf{Bottom:} Structure solver.}
    \label{fig:sech}
\end{figure}

As shown in Fig. \ref{fig:sech}, in the fluid dynamic solver, we are interested in learning the operator mapping the velocity from $t_0$ up to time $t_0+\alpha$ to the velocity up to some later time $T>t_0+\alpha$, $G^{\dagger} :C([t_0,t_0+\alpha]; H_{per}^r((-1,1)^3;\mathbb{R}))) \to C((t_0+\alpha,T];H_{per}^r((-1,1)^3;\mathbb{R}))$ defined by ${\bf u}|_{(-1,1)^3\times [t_0,t_0+\alpha]}\to$ ${\bf u}|_{(-1,1)^3\times (t_0+\alpha,T]}$. The operator maps the solution at the previous $\alpha$ time steps to the next time step. For example, the data from the dashed yellow box is used to predict ${\bf u}^{t_0+\alpha+2}$, with a yellow arrow connecting the feature and label. The FNO is achieved using a combination of 3D-convolutional layers for spatial information processing and recurrent layers for temporal information propagation. So the FNO solver obtains the flow field at $t_0 +\alpha +1 \sim T$. Then the velocity of fluid nodes is transferred to the vesicle as per Eq. (\ref{coupling}). Only needing the position of structure at ${\bf X}_k^{t_0+\alpha}$, we get the position at $t_0+\alpha+1 \sim T$ by iteratively solving Eq. (\ref{coupling}) using a finite different method,
\begin{equation}\label{vesicleSim}
{\bf X}_k^{n+1}={\bf X}_k^{n}+\sum_{\bf m}\delta_a({\bf x}_{\bf m}-{\bf X}_k^{n}) {\bf u}_{\bf m}^{{\rm ind},n} \Delta x^3 + {\bf u}^{\infty}({\bf X}_k^{n}) \Delta t, 
\end{equation}
where the detailed Dirac-delta function is seen from \cite{PhysRevE.107.035103}. The fluid-structure interaction problem for the time interval $t_0+\alpha+1$ to $T$ is fully solved.

\section{Results and discussions}
\subsection{Data preparation}
As mentioned in Fig. \ref{fig:sech}, only the fluid flow is predicted by the FNO. The structure is modeled by an integral equation and solved by a differential scheme. We train the FNO for interaction with a vesicle on simulations with a prescribed motion of the vesicle with two different initial shapes $\{{\bf X}^0_k\}^{M}_{k=1}$.

For each initial vesicle $\{{\bf X}^0_k\}^{M}_{k=1}$, the equations (\ref{fluidv}) - (\ref{coupling}) are solved using the IBM presented by Xiao {\it et al.} \cite{PhysRevE.107.035103}. All data are generated on a $64^3$ grid. The shear rate of the elongation flow is 1. We use a time-step of $\Delta t = 10^{-6}s$ in the data-generated process where we record the solution every $\Delta T = 2\times 10^{-3}s$ time units. After $0.5s$ of turning on the elongation flow, we observe that the vesicle reaches an equilibrium state. We control the elongation flow suddenly switched at $N\Delta T = 0.604s$ with $N=302$. Thus we study the process after flow reversal in this paper. We observe that the system reaches the new steady state around $N\Delta T = 1.04s$ and $N=502$. We store $1000$ time steps of the data for each simulation. We now preprocess the collected trajectories into the needed format for training.

For time series prediction, we choose the input format to be a sliced window $\alpha = 10$, as shown in Fig. \ref{fig:sech}. We consider the label in two cases: one-step label and multi-step labels. We investigate the effectiveness of these two forms, referred to as 'sequence to one' and 'sequence to sequence' respectively. The training data is also divided into two cases: flow velocity data without steady state (specifically from $N=302$ to $N=481$) and with steady state (specifically from $N=302$ to $N=521$). Thus we establish four distinct experiments, as illustrated in Table \ref{tab:table1}, as follows:
% 本文中，我们考虑用10个时刻的流速作输入数据去预测接下来1个时刻和接下来10个时刻的流速，研究这两种预测形式的效果。为了简化描述形式，分别称这两种预测形式为'10 to 1' 和 '10 to 10'形式。训练数据也分成两种，一种是取到未达到稳态的流速数据，即第302-481时刻的流速；另一种是取到达到稳态的流速数据，即第302-521时刻的流速。研究两种样本量的预测效果。基于这两种考虑，我们设置了如表一所示的四种网络。下面给出了具体表述：

\begin{itemize}
    \item FNO-1: The training data includes flow velocity \textbf{without} steady state which is trained in the '\textbf{sequence to one}' form. We have $170$ samples for each vesicle and a total of $340$ training samples for two vesicles.

    \item FNO-2: The training data includes flow velocity \textbf{with} steady state which is trained in the '\textbf{sequence to one}' form. We have $210$ samples for each vesicle and a total of $420$ training samples.

    \item FNO-3: The training data includes flow velocity \textbf{without} steady state which is trained in the '\textbf{sequence to sequence}' form.
    % For example, we sequentially select 20 time steps of flow velocity starting from $N=302$. The flow velocity at $N=302-311$ is taken as the input, while the flow velocity at $N=312-321$ is taken as the output, forming the first sample data of the training dataset. 
    We have $160$ samples for each vesicle and a total of $320$ training samples.

    \item FNO-4: The training data includes flow velocity \textbf{with} steady state which is trained in the '\textbf{sequence to sequence}' form. We have $200$ samples for each vesicle and have a total of $400$ training samples.
\end{itemize}

\begin{table}[t]
\begin{center}
\caption{\label{tab:table1}Four experiments of operator.}
% \begin{ruledtabular}
\resizebox{1\columnwidth}{!}{
\begin{tabular}{ccccccc}
\multirow{2}{*}{Experiment}&\multirow{2}{*}{Train method} & \multirow{2}{*}{Train Dataset}&\multicolumn{1}{c}{Training }&\multicolumn{3}{c}{Prediction Type} \\
      &                      &         & Sample & inter-set & mix-set & extra-set \\ \hline
 FNO-1& one-step label     & without steady state & 340 & \multirow{4}{*}{interpolation} & extrapolation & \multirow{4}{*}{extrapolation} \\ 
 FNO-2& one-step label     & with steady state & 420 &  & interpolation &  \\
 FNO-3& multi-step labels & without steady state & 320 &  & extrapolation & \\
 FNO-4& multi-step labels & with steady state & 400 &  & interpolation &  \\ \hline
\end{tabular}
}
% \end{ruledtabular}
\end{center}
\end{table}

We construct the FNO by stacking four Fourier integral operator layers as specified in (\ref{FNO1}) and (\ref{FNO3}) with the ReLU activation. We use the Adam optimizer to train for 200 epochs with a learning rate of 0.001. We set $k_{max} = 12$, $d_v=64$. All the computation is carried on a single Nvidia A100 GPU with 80GB memory.
% 对于FNO-1，我们以如下形式重组数据构建训练集。对每个仿真过程，从$N=302$个时刻开始按顺序取11个时刻的流速解，第302-311时刻的流速为输入，第312时刻的流速作为输出，构成训练集的第1个样本数据。然后从$N=303$个时刻开始按顺序取11个时刻的解，第303-312个时刻的流速作为输入，第313时刻的流速作为输出，构成训练集的第2个样本数据。依次类推，这个过程直到取到第471个时刻，共取170组数据作为训练集。两个vesicle就会取得共320组数据作为训练集来训练FNO-1。对于FNO-2，训练集的构成和FNO-1的类似，只是会一直取到$N=521$时刻的流速共420组样本数据。FNO-3的训练集的选取与FNO-1的类似。只是训练集的样本由二十个时刻的流速构成。前十个时刻的流速作输入，后十个时刻的流速是输出。共有320组样本数据。 FNO-4的训练集的选取与FNO-3的类似，一共有400组样本数据。

% 对于FNO-3，我们以如下形式重组数据构建训练集。对每个仿真过程，从$N=302$个时刻开始按顺序取20个时刻的流速解，第$N=302-311$时刻的流速为输入，第$N=312-321$时刻的流速作为输出，构成训练集的第1个样本数据。然后从$N=303$个时刻开始按顺序取20个时刻的解，第$N=303-312$个时刻的流速作为输入，第$N=313-322$时刻的流速作为输出，构成训练集的第2个样本数据。依次类推，取到第$N=461$个时刻，共取160组数据作为训练集。两个vesicle就会取得共320组数据作为训练集来训练FNO-3。对于FNO-4，训练集的构成和FNO-3的类似，只是会一直取到$N=521$时刻的流速共400组样本数据。

\subsection{Performance comparisons on flow prediction}
\begin{figure}[h]
     \centering    
     \includegraphics[width=1\linewidth]{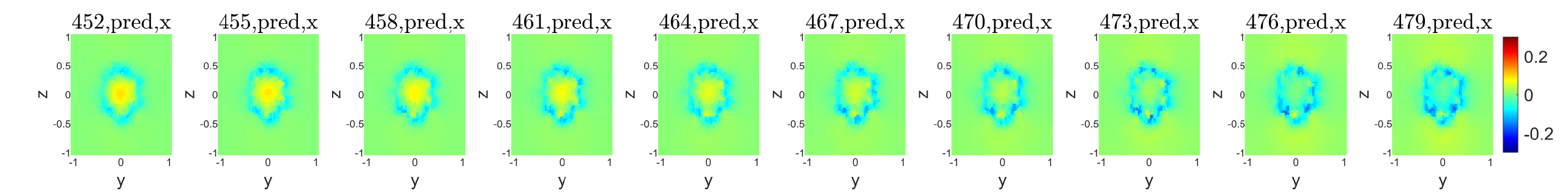}(a)\\
     \includegraphics[width=1\linewidth]{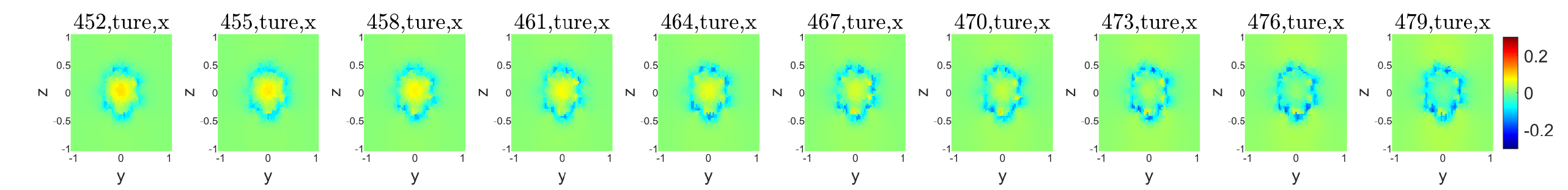}(b)\\
     \includegraphics[width=1\linewidth]{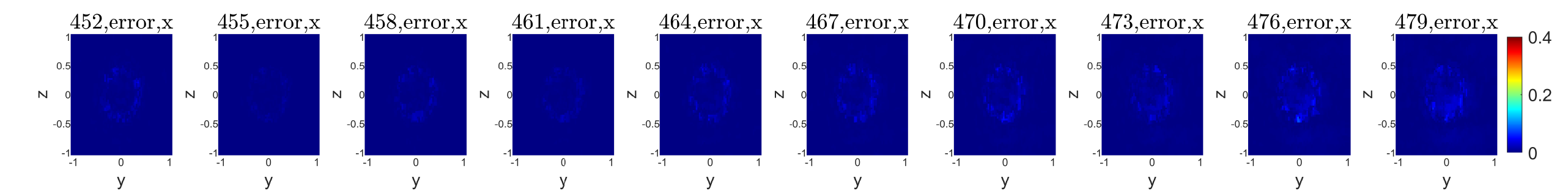}(c)\\
     \caption{(a) the predicted results and (c) errors of the FNO-4 with inter-set as input. Here, the error refers to the absolute difference between the predicted solution and the IBM solution (b).}    
     \label{exp_fluid} 
 \end{figure}

We test the performance of the four experiments mentioned above, using the velocity data at ten-time steps $442$-$451$, $482$-$491$, and $522$-$531$ as input. For convenience, these three inputs are named inter-set, mix-set, and extra-set, respectively. Figure \ref{exp_fluid} shows the predicted results and errors of FNO-4 with inter-set as input. We only exhibit the x-component of the velocity. The error refers to the absolute difference between the predicted solution of the FNO and the IBM solution. The other predicted results and errors of the four sets mentioned above are shown in Appendix. The experiments with inter-set and extra-set as input are interpolation and extrapolation, respectively, as shown in Table \ref{tab:table1}. Other experiments with mix-set as input are interpolation or extrapolation.

% Figure shows the average errors of the predicted solutions for the 30 time steps across all FNOs.
% 以$N=442-451$，$N=482-491$，和$N=522-531$十个时刻的流速做输入数据，上述4个网络的预测结果及误差分别如图\ref{exp442_451}，\ref{exp482_491}和\ref{exp522_531}所示。这些是x方向上流速的演示图。这里的误差是指预测解与IBM解差值的绝对值。每张图的子图中，(a)是IBM求得的解。(b), (d), (f), and (l)分别是模型1，模型2，模型3，和模型4的预测解。(c), (e), (h), and (k)分别是图例最大值为0.4的FNO-1，FNO-2，FNO-3，和FNO-4的误差。(g) and (m)分别是图例最大值为0.05的模型3和模型4的误差。图\ref{fluid_mean_err}展示的是各个模型预测的30个时刻的解的平均误差。

\subsubsection{'Sequence to one' vs 'sequence to sequence'}
\begin{figure}[h]
    \centering
    \includegraphics[width=0.5\linewidth]{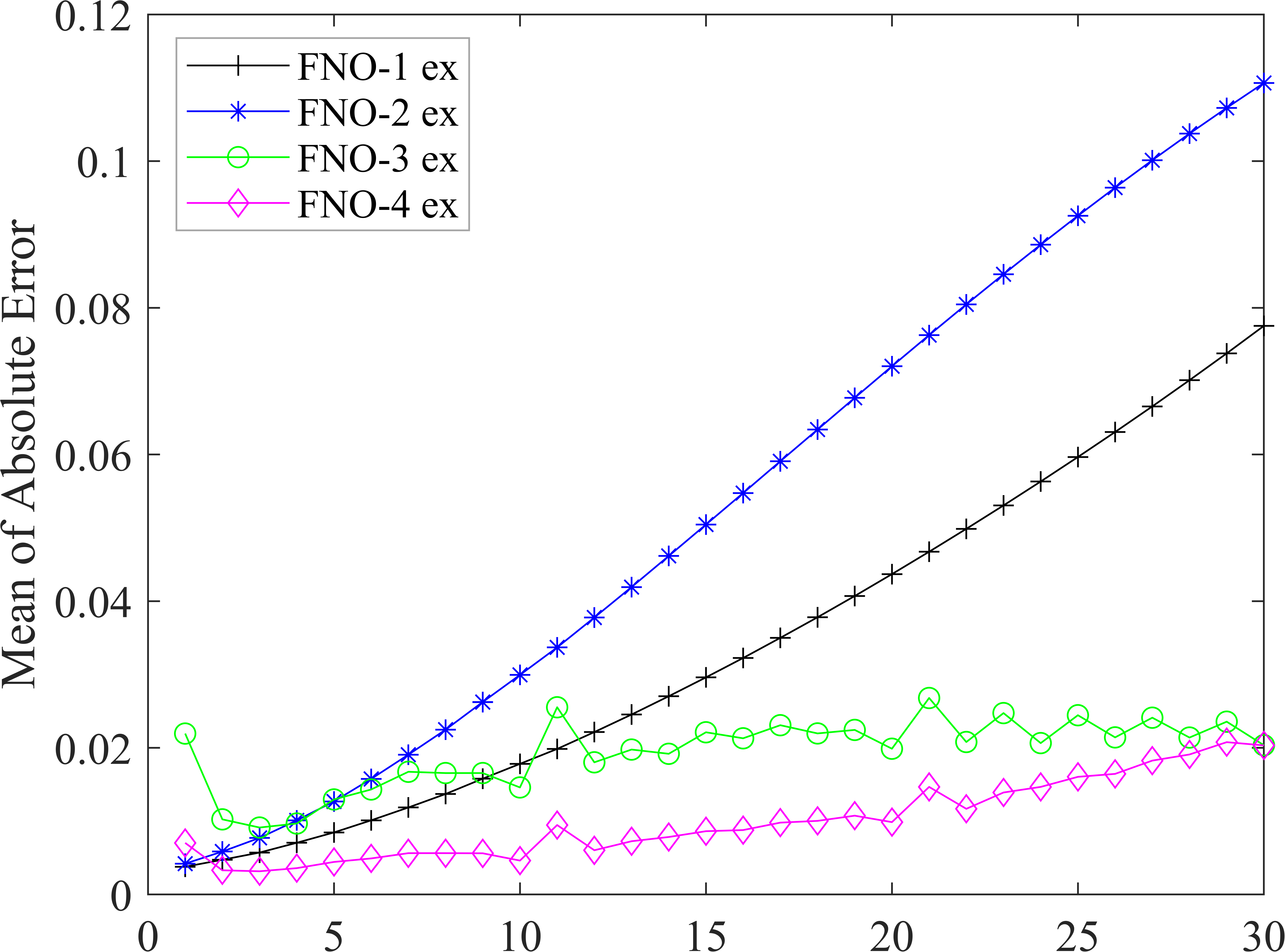}
    \caption{The mean of absolute error between the predicted solution and the IBM solution for the predicted 30-time steps is shown. This figure only demonstrates the x-component of the velocity. The 'ex' in the legend represents the prediction error with extra-set as input.}
    \label{one_vs_seq}
\end{figure}

We investigate the effects of 'sequence to one' and 'sequence to sequence' training on FNO prediction. Figure \ref{one_vs_seq} shows the error of four experiments with extra-set as input. These predictions of the flow velocity are all results of extrapolation. The errors of FNO-1 and FNO-2 (sequence to one) are smaller than those of FNO-3 and FNO-4 (sequence to sequence) at the first time step of prediction. Their errors become larger in subsequent time steps. These phenomena indicate that the 'sequence to one' prediction format is more suitable for short-term predictions. The error curves of FNO-1 and FNO-2 (sequence to one) are smooth and increasing which indicates that their errors strictly increase with time. On the other hand, FNO-3 and FNO-4 (sequence to sequence) do not exhibit such a pattern, as their error curves fluctuate and increase.

The long-term prediction capabilities are examined by the accuracy of the predicted flow velocities from the 15th to the 30th time steps. Figure \ref{one_vs_seq} shows that FNO-3 and FNO-4 (sequence to sequence) have higher accuracy than FNO-1 and FNO-2 (sequence to one). The experiments with inter-set as input have similar results. 

The above results show that the 'sequence to sequence' type FNO performs better than the 'sequence to one' type FNO for long-term predictions. This result is because the errors of the 'sequence to one' type FNO accumulate over time and eventually become consistently larger than the errors of the 'sequence to sequence' type FNO. In the rest of this paper, we only discuss the 'sequence to sequence' type FNO. 

% FNO-3 and FNO-4 have higher accuracy than FNO-1 and FNO-2, as shown in Figs. \ref{exp442_451} and \ref{fluid_mean_err}. The predictions with extra-set are all results of extrapolation. Figures \ref{exp522_531} and \ref{fluid_mean_err} illustrate that FNO-4 and FNO-3 also have higher accuracy than FNO-1 and FNO-2. When predicting in mix-set, Figs. \ref{exp482_491} and \ref{fluid_mean_err} demonstrate that FNO-4 have higher accuracy than FNO-2, which are both interpolations. FNO-3 have higher accuracy than FNO-1, which are both extrapolations. 

% 从图\ref{fluid_mean_err}中可以看出来，第一个时刻，FNO-1和FNO-2的误差比FNO-3和FNO-4的小，但其后的误差都要更大。说明'sequence to one'预测形式的更适用做短时预测。在图\ref{fluid_mean_err}中，FNO-1和FNO-2的误差曲线是光滑上涨的。也就是说其误差是严格随时间增大的。而FNO-3和FNO-4却没有类似的规律。他们的误差曲线是波动上涨的。选取图\ref{fluid_mean_err}中25时刻后的数据，横向对比输入数据相同的曲线发现，不论训练集是否取到了稳态时刻的流速，也不论是做内插还是外插，'sequence to sequence'型FNO的误差比'sequence to one'型FNO小。所以说，从长时预测的角度来说，'sequence to sequence'型FNO的预测效果会比'sequence to one'的型FNO的效果更好。这是因为'sequence to one'型FNO的误差严格随时间累计，累计到一定程度就会一直大于'sequence to sequence'型FNO的误差.这一点从图\ref{exp442_451}-\ref{exp522_531}中也可以看出来。

% Regarding the solutions after the 24th time step, the accuracy of models trained in the 'sequence to one' form for interpolation is still higher than that of models trained in the 'sequence to sequence' form for extrapolation, but it is approaching similar levels, and even smaller between the 13th and 18th time steps.

\subsubsection{Interpolation vs extrapolation}
\begin{figure}[h]
    \centering
    \includegraphics[width=0.45\linewidth]{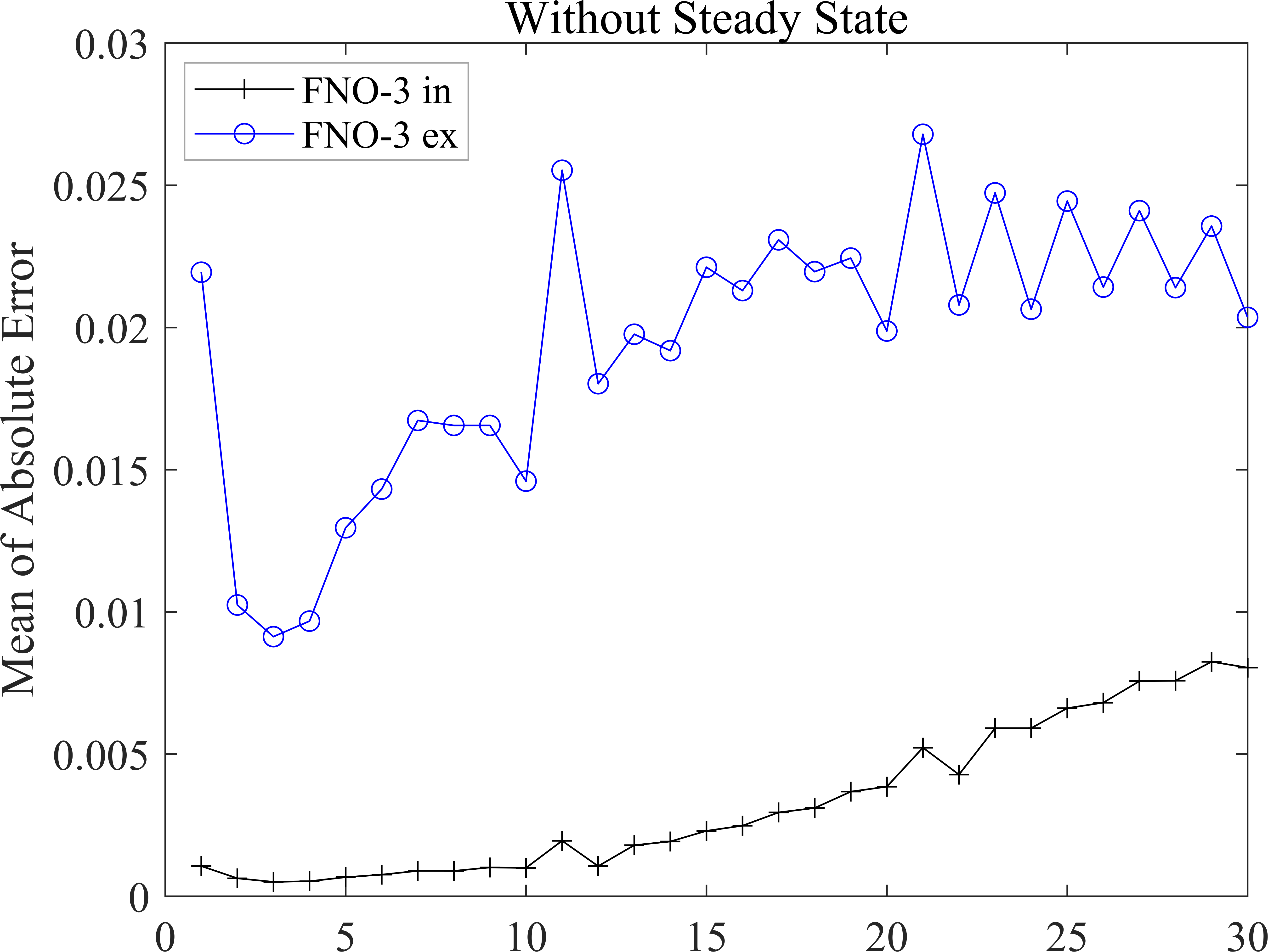} (a)
    \includegraphics[width=0.45\linewidth]{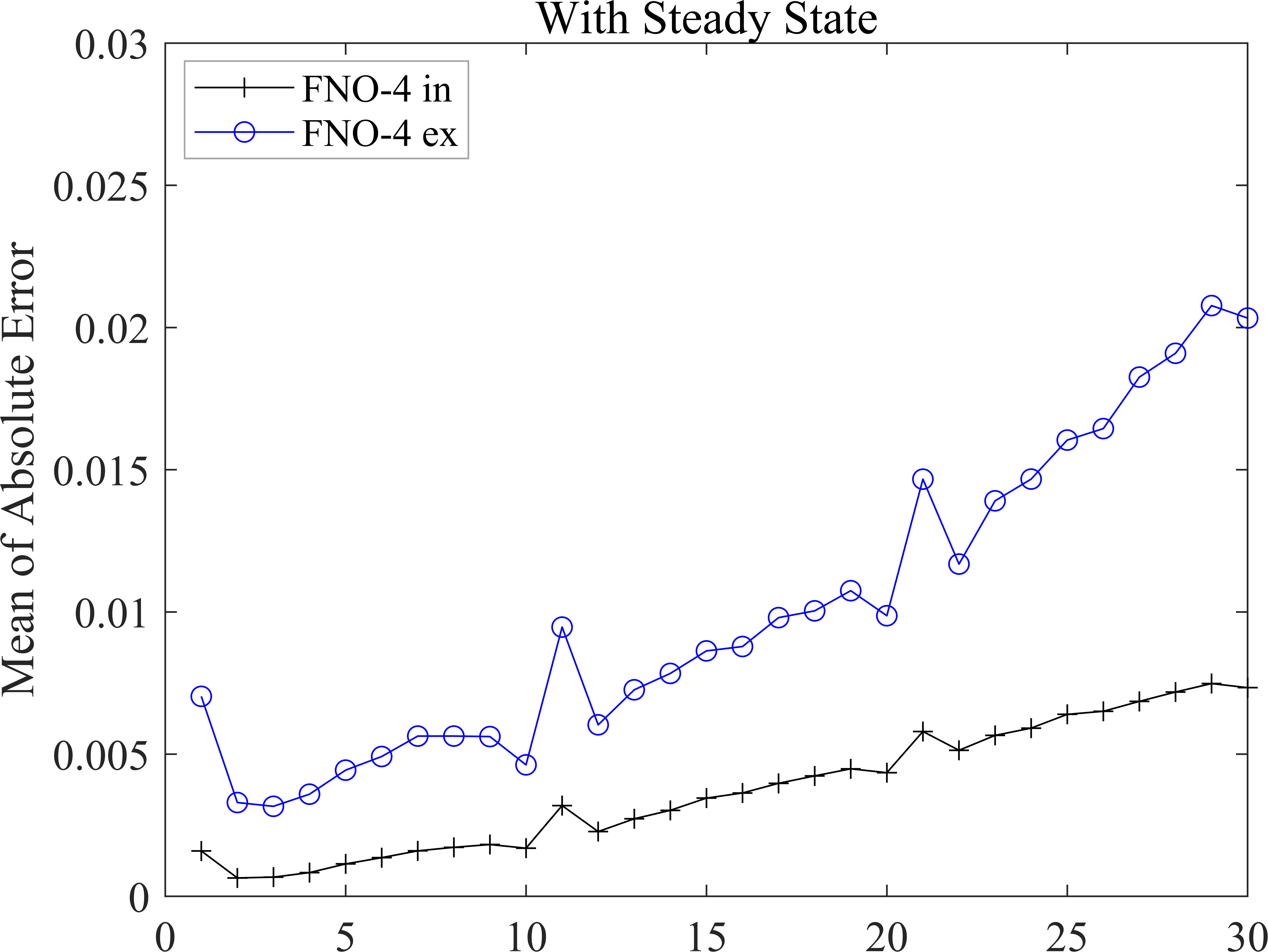} (b)
    \caption{The mean of absolute error of experiments (a) without steady state and (b) with steady state. The 'in' in the legend represents the error of the prediction experiment with inter-set as input.}
    \label{in_vs_ex}
\end{figure}

We analyze the prediction accuracy of interpolation and extrapolation of 'sequence to sequence' type FNO on the same dataset. The experiments with inter-set and extra-set as input are interpolations and extrapolation, respectively. Figure \ref{in_vs_ex} shows that the performance of interpolation is better than extrapolation, both without (Fig. \ref{in_vs_ex}(a)) and with (Fig. \ref{in_vs_ex}(b)) steady-state conditions.  

% The training datasets of FNO-4 include the steady-state velocity. Their prediction errors are minimized when using the inter-set as input. The predictions error is somewhat higher when using the inter-set as input. It reaches its maximum when employing the extra-set as input. On the other hand, the training datasets of FNO-1 and FNO-3 do not include the steady-state velocity. Their prediction errors are minimized when using that at $N=442-451$ as input, which are both interpolations. These results demonstrate that the performance of interpolation is better than extrapolation regardless of the scenario.

% The prediction error of FNO-1 and FNO-3 when using the inter-set as input is greater than that when using the extra-set as input, which are both extrapolations. This maybe due to that the fluid undergoes larger changes at $N=492$−$521$ compared to the already steady-state changes at $N=532$−$561$. Therefore, the accuracy of predicting the former is lower than the latter.

%不论是哪种情况，做内插的效果都要优于做外插的效果。若将训练数据集取到稳态时刻的流速，即看FNO-2和FNO-4的结果，发现以$N=442-451$时刻流速作为输入的预测误差最小。以$N=482-491$时刻流速作为输入的预测误差次之。以$N=522-531$时刻流速作为输入的预测误差最大。而若训练集未取到稳态时刻，那么预测$N=442-451$时刻的误差仍是最小。但是以$N=482-491$时刻流速为输入的的预测误差比以$N=522-531$时刻流速为输入的更大。这可能是因为$N=482-491$时流体变化比已经达到稳态的$N=522-531$时的变化大，所以同样是做外插时，前者的准确率反而比后者的低。

\subsubsection{With steady states vs without steady states}
\begin{figure}[h]
    \centering
    \includegraphics[width=0.45\linewidth]{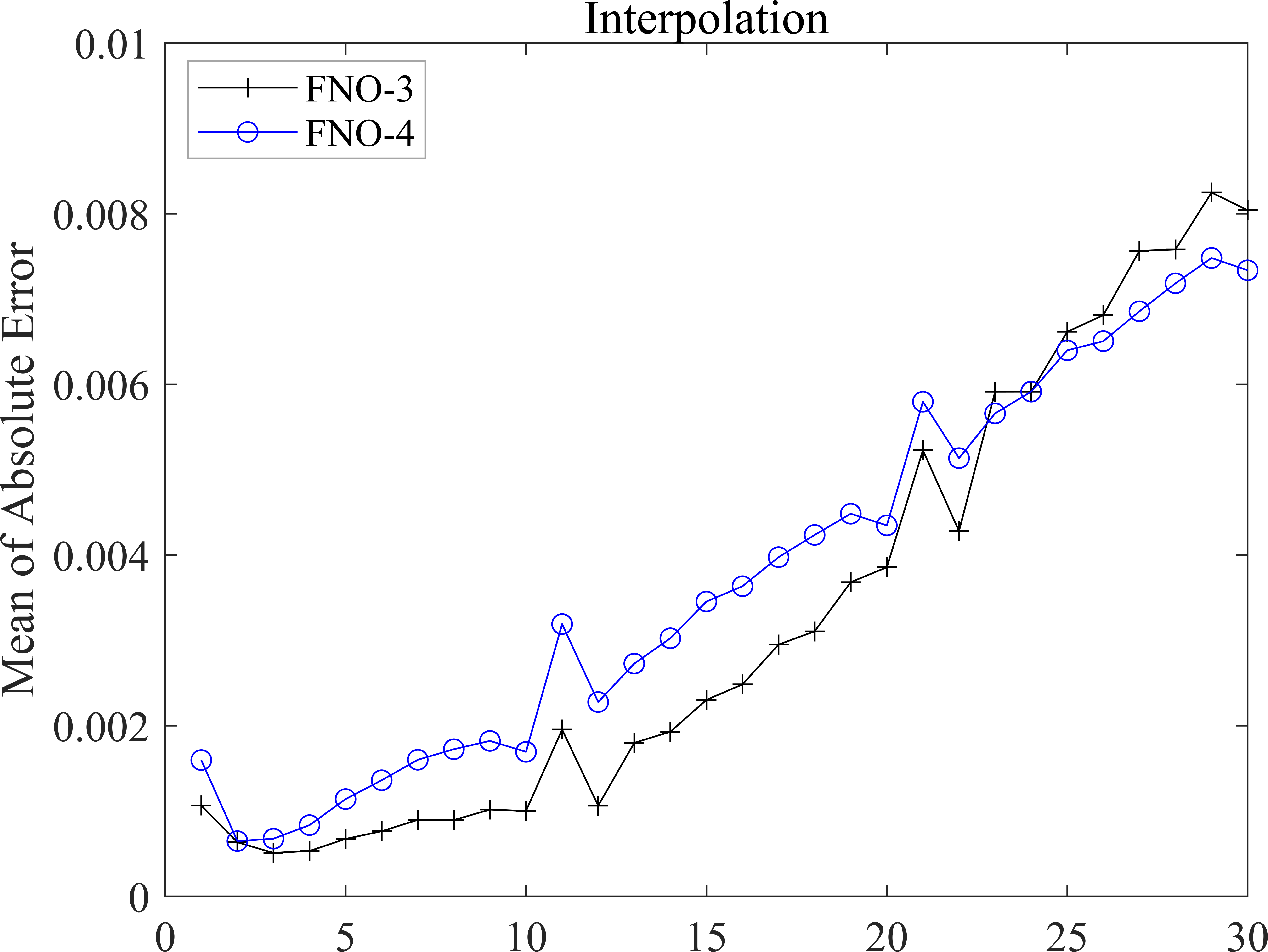} (a)
    \includegraphics[width=0.45\linewidth]{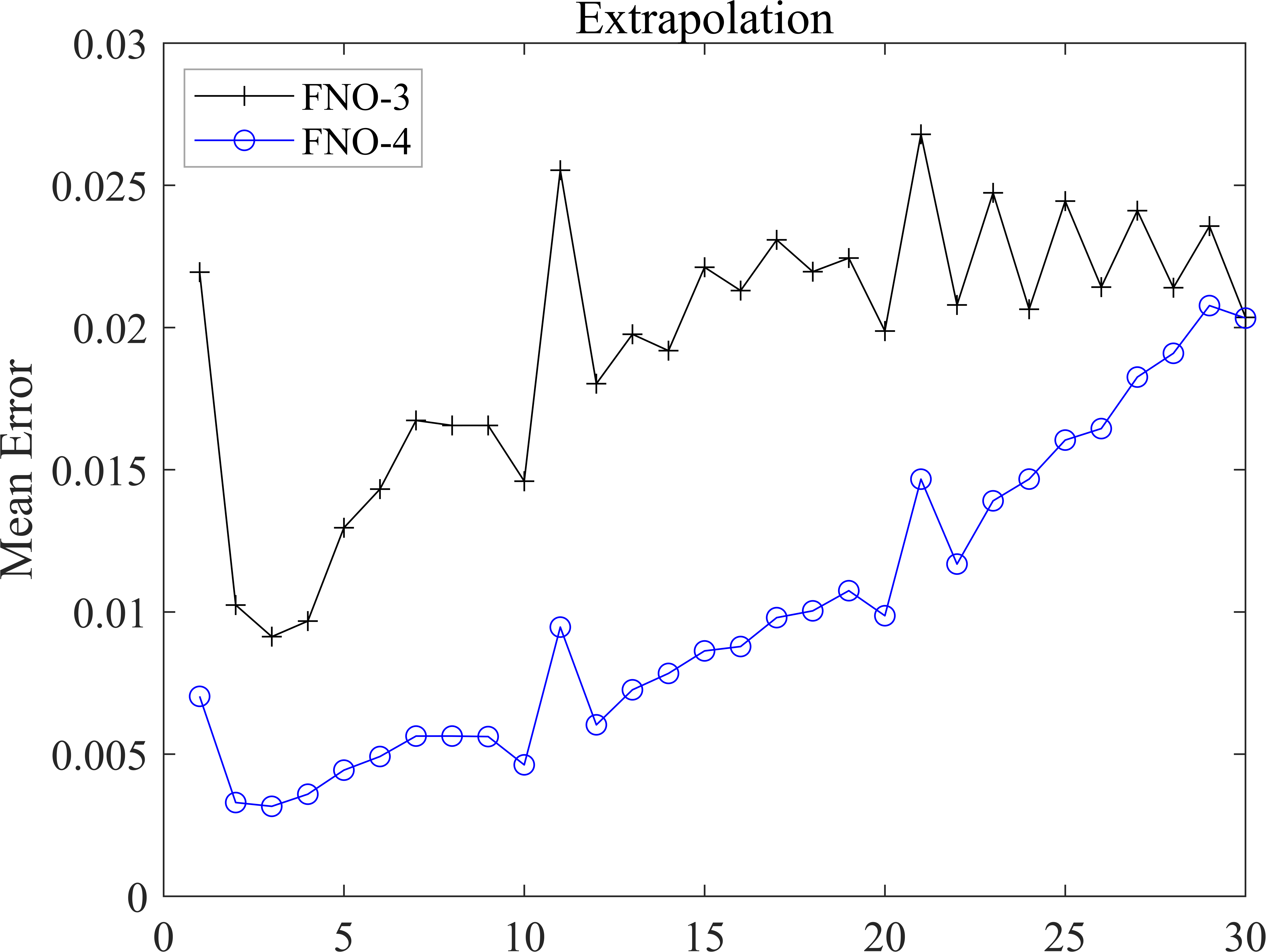} (b)
    \caption{The mean of absolute error of experiments of (a) interpolation and (b) extrapolation.}
    \label{with_vs_without}
\end{figure}

We analyze the prediction accuracy of different models for the different training sets. Figure \ref{with_vs_without} demonstrates that when predicting the flow velocity for the subsequent 25-time steps, FNO-4 (with steady states) performs better than FNO-3 (without steady states), both for interpolation (Fig. \ref{with_vs_without}(a)) and for extrapolation (Fig. \ref{with_vs_without}(b)). These results indicate that higher accuracy is obtained when the training set without steady-state data.

% there are all results of interpolation. FNO-1 performs better than FNO-2, and FNO-3 performs better than FNO-4. These results indicate that higher accuracy is obtained when the training set without steady-state data. Therefore, we conclude that the more accurate the training set, the higher the interpolation accuracy. {\color{red} From the 25th time step onwards, FNO-4 outperforms FNO-3 in terms of accuracy.}

% When predicting that after the extra-set, there are all results of extrapolations. Figures \ref{exp522_531} and \ref{fluid_mean_err} illustrate that FNO-1 performs better than FNO-2, but FNO-4 performs better than FNO-3. No consistent conclusion can be drawn.

% The long-term prediction capabilities are examined by the accuracy of the predicted flow velocities for the 15th to the 30th time steps.

% 我们接下来在相同数据集的基础上，分析不同模型的预测精度。图1和图4表明，用四个算子预测$N=442-451$后22个时刻的流速，精度最高的是FNO-3，其次是FNO-4，然后是FNO-1，最差的是FNO-2。FNO-1的效果比FNO-2的好，FNO-3的效果比FNO-4的好，这说明训练集不取到稳态的预测精度反而更高。而这四组数据都是做内插的结果。那么我们可以说训练集越精准，内插精度越高。在第25个时刻及以后，FNO-4的精度优于FNO-3。对于$N=522-531$时刻的预测结果，我们发现，精度最高的是FNO-4，其次是FNO-3，然后是FNO-1，最差的是FNO-2。FNO-1的效果比FNO-2的好，但FNO-4的效果比FNO-3的好.得不出有规律的结论。对于482-491时刻的预测结果，同样是内插的FNO-2和FNO-4的精度分别比同样是外插的FNO-1和FNO-3的精度高。针对24个时刻及之后的解来说，以'sequence to one'形式训练的模型做内插的精度虽然仍是大于以'sequence to sequence'形式训练的模型做外插的精度，但已经接近了，在13-18个时刻间甚至还更小。

\subsection{Vesicle}
\begin{figure}[h]
     \centering
     \includegraphics[width=1\linewidth]{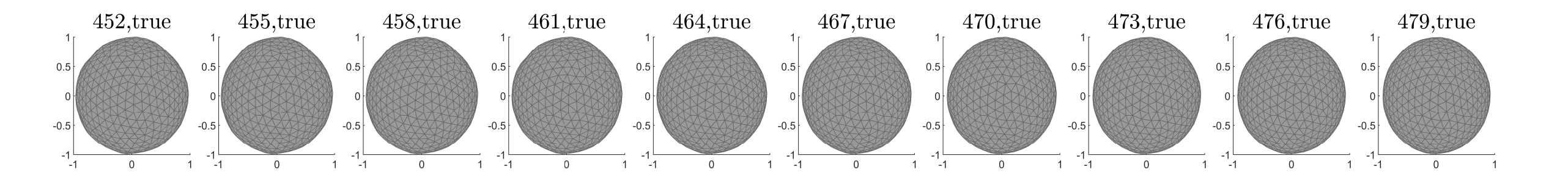}(a)\\
     \includegraphics[width=1\linewidth]
     {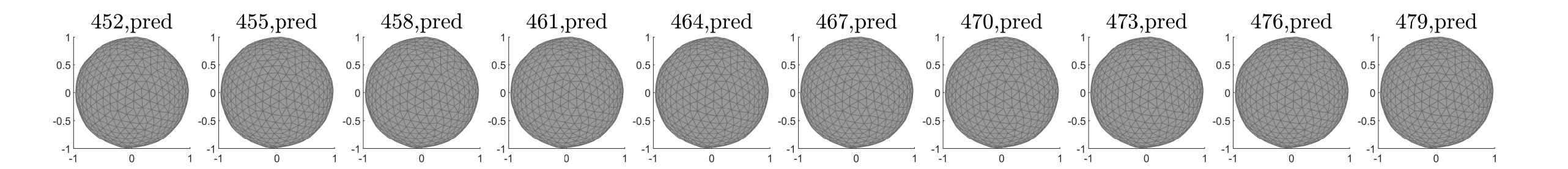}(b)\\
     \includegraphics[width=1\linewidth]{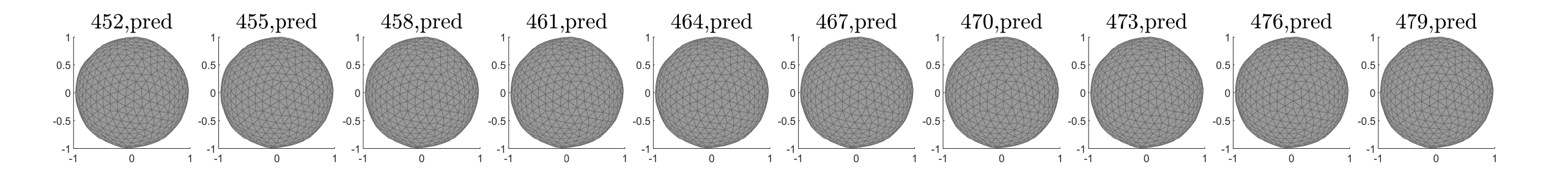}(c)
     \caption{The evolution of the vesicle at the last 30 moments with inter-set as input. (a) represents the solution obtained by IBM. (b) represents the solution obtained by using the fluid velocity predicted by FNO-3. (c) represents the solution obtained by using the predictions by FNO-4.}
     \label{v_452}
 \end{figure}

The 'sequence to sequence' approach has a higher accuracy for long-term prediction. 
% As shown in Fig. \ref{exp_fluid}(a), we plot the predictions of FNO-4 with inter-set as input. Figure \ref{exp_fluid}(b) show the results of the IBM. Figure \ref{exp_fluid}(c) shown the absolute difference between the predicted solution and the IBM solution. 
We transfer these predictions of FNO to Eq. (\ref{coupling}) to solve the dynamic of the vesicle by Eq. (\ref{vesicleSim}). Figure \ref{v_452} shows the evolution of the vesicle at the last 30 moments with inter-set as input. Here, Fig. \ref{v_452}(a) represents the solution obtained by transferring the fluid velocity predicted by the IBM. Figure \ref{v_452}(b) represents the solution obtained by using the predictions by FNO-3. And Fig. \ref{v_452}(c) shows the results of Eq. (\ref{vesicleSim}) based on Fig. \ref{exp_fluid}(b). The evolution of the vesicle at the last 30 moments with mix-set and extra-set as input are shown in Fig. \ref{v_492} and Fig. \ref{v_532}, respectively, in the Appendix.

% The evolution of vesicles at $N=452$−$481$, $N=492$−$521$, and $N=532$−$561$ are shown in Figs. \ref{v_452}, \ref{v_492}, and \ref{v_532}, respectively. Subfigure (a) represents the solution obtained by transferring the fluid velocity predicted by the IBM. Subfigure (b) represents the solution obtained by transferring the prediction of FNO-3, and subfigure (c) represents the solution obtained by transferring the prediction of FNO-4.

\begin{figure}[h]
     \centering
     \includegraphics[width=0.5\linewidth]
     {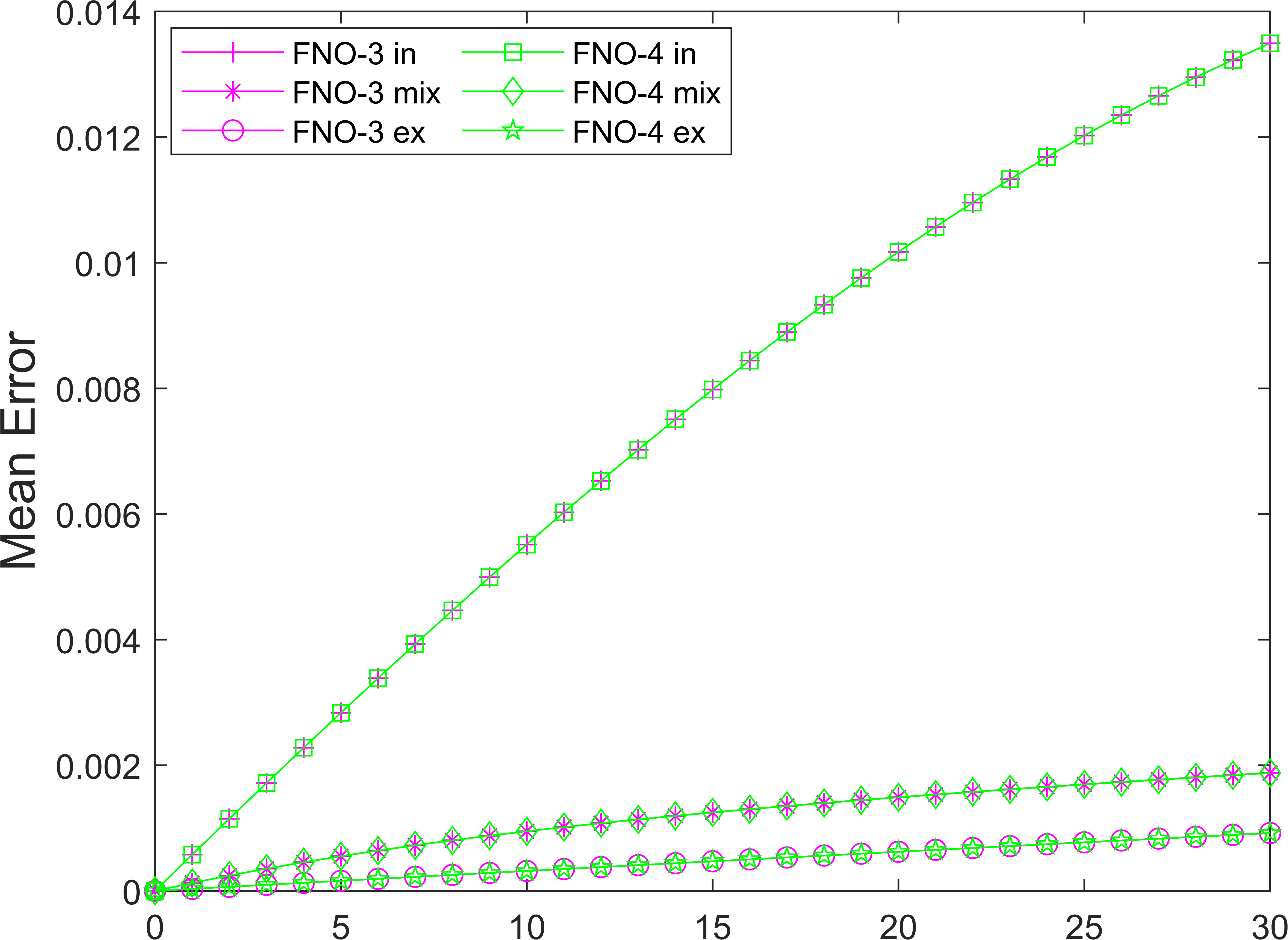}
     \caption{Graph of the variation of the mean error of the vesicles over time. The 'mix' in the legend represents the error of the prediction experiment with mix-set as input.}
     \label{vesicle_mean_err}
 \end{figure}

Figure \ref{vesicle_mean_err} shows the mean errors of the aforementioned results. The different training sets have almost no impact on the accuracy of vesicle prediction. The highest prediction accuracy is achieved at the last 30 moments with extra-set as input. This result indicates that the model is most accurate in predicting the steady state. The lowest prediction accuracy is obtained at the last 30 moments with inter-set as input. It is worth noting that the fluid prediction accuracy for this interval is the highest. This may be due to the greater influence of vesicle deformation on the error in this interval compared to the other two intervals.

\section{Conclusion}
In this paper, we introduce an FSI solver that uses the FNO to predict fluid flow and apply the FSI solver to the problem of interaction between elongation flow and vesicle. The model is trained in 3D space with both the position of fluid flow and vesicle, as well as the velocity of the fluid.

Additionally, we investigate the FNO's ability to multi-step prediction, so compare two cases: 'sequence to one' and 'sequence to sequence'. In order to validate the performance of the FNO-based FSI solver, we conduct a couple of simulations for these two scenarios. We compared the predictions with the results of the IBM. We found that 'sequence to one' types of FNO achieve higher accuracy in short-term prediction, while 'sequence to sequence' types of FNO demonstrated higher accuracy in long-term prediction. The results obtained by the IBM and 'sequence to sequence' types of FNO-based FSI solvers exhibit similar quantitative behavior. Specifically, the fluid variations near the surface of the membrane are larger and eventually reach an equilibrium state. We also investigated the training effect of datasets that include steady-state solutions and those that do not include steady-state solutions. Moreover, we find that the performance of interpolation is better than extrapolation regardless of the scenario.

Finally, we utilized the predictions from 'sequence to sequence' types of FNO and simulated the dynamic process of the vesicle using a solid solver. For the solid solver, we modeled the dynamics of the structure by an integral equation equation using IBM. This method utilizes the Dirac-delta function to transfer the velocity from fluid to the structure. We obtained the position of the structure by solving the integral equation using a finite different method.

The proposed FNO-based FSI model shows good potential in structural optimization or flow control in problems where FSI effects are evident. However, its prediction accuracy and generalization still need to be further improved. In the next step, we will further optimize the FNO-based FSI model to obtain a more accurate interaction between the fluid and the structure. Embedding the Navier–Stokes equation and structural motion control equations into the neural improves its generalization ability. The new FNO-based FSI model is expected to be used in the FSI problems with different incoming flow states. Moreover, it is expected to simulate the wrinkling of a vesicle.

\section*{CRediT authorship contribution statement}
\textbf{Wang Xiao:} Conceptualization, Methodology, Software, Formal analysis, Investigation, Writing - Original Draft. 
\textbf{Ting Gao:} Conceptualization, Methodology, Software, Formal analysis, Writing - Review \& Editing, Resources.
\textbf{Kai Liu:} Resources, Writing - Review \& Editing, Supervision, Funding acquisition. 
\textbf{Meng Zhao:} Methodology, Validation, Formal analysis, Resources, Writing - Review \& Editing, Supervision, Funding acquisition.
\textbf{Jinqiao Duan:} Project administration, Resources, Supervision, Funding acquisition.

\section*{Declaration of competing interest}
The authors declare that they have no known competing financial interests or personal relationships that could have appeared to influence the work reported in this paper.

\section*{Data availability}
Data will be made available on request.

\section*{Acknowledgments}
This work was supported by the National Key Research and Development Program of China (No. 2021ZD0201300), the National Natural Science Foundation of China (No. 12141107 and No. 12301553), the Fundamental Research Funds for the Central Universities, China (5003011053) and Dongguan Key Laboratory for Data Science and Intelligent Medicine, China. M. Z. also thanks the support from the Ministry of Education Key Lab in Scientific and Engineering Computing.

\appendix
\section{The predicted results of four sets}
\begin{figure}[h]
     \centering    
     \includegraphics[width=1\linewidth]{stead/steady_state_in_452_ture.png}(a)\\
     \includegraphics[width=1\linewidth]{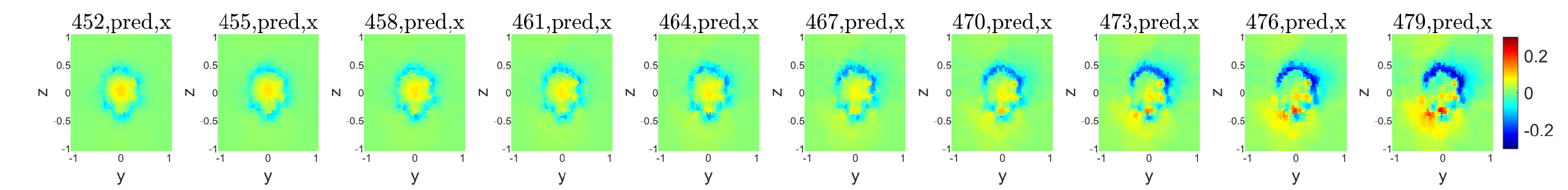}(b)\\
     \includegraphics[width=1\linewidth]{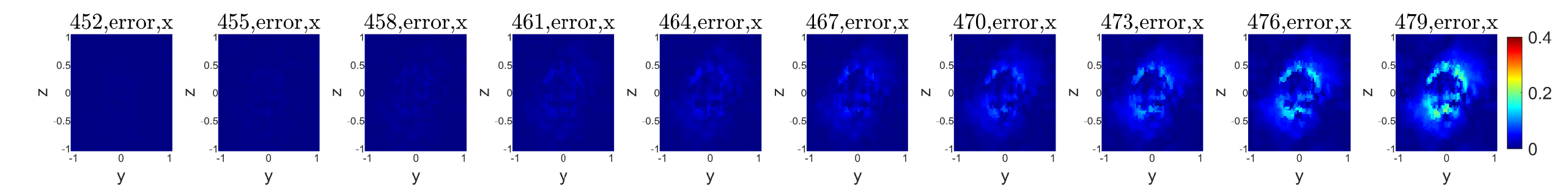}(c)\\
     \includegraphics[width=1\linewidth]{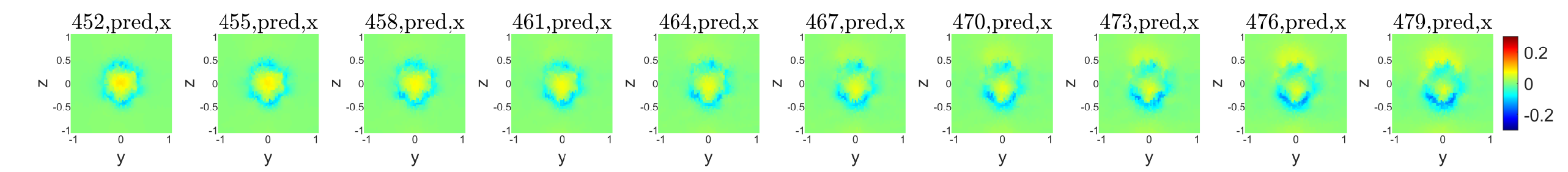}(d) \\
     \includegraphics[width=1\linewidth]{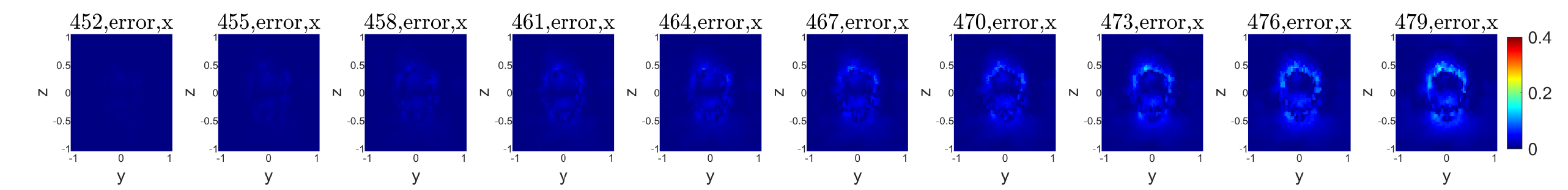}(e)\\
     \includegraphics[width=1\linewidth]{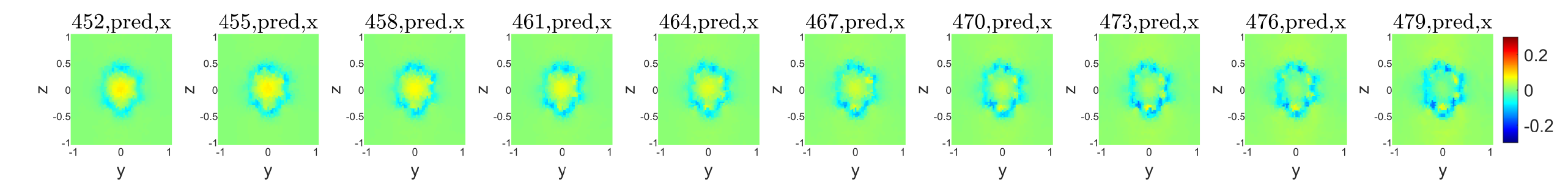}(f)\\
     \includegraphics[width=1\linewidth]{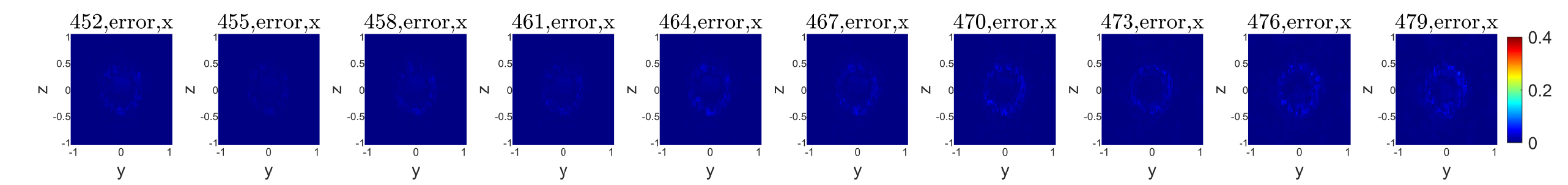}(g)\\
     \includegraphics[width=1\linewidth]{10to10/to10_steady_in_452_pred.png}(h)\\
     \includegraphics[width=1\linewidth]{10to10/to10_steady_in_452_err_4.png}(l)\\
     \caption{The predicted results and errors of the four FNOs with inter-set as input are shown. (a) represents the solution obtained by IBM. (b), (d), (f), and (h) represent the predicted solutions of FNO-1, FNO-2, FNO-3, and FNO-4, respectively. (c), (e), (g), and (l) represent the errors of FNO-1, FNO-2, FNO-3, and FNO-4, respectively. Here, the error refers to the absolute difference between the predicted solution and the IBM solution.}     
     \label{exp442_451} 
 \end{figure}

\begin{figure}[h]
     \centering    
     \includegraphics[width=1\linewidth]{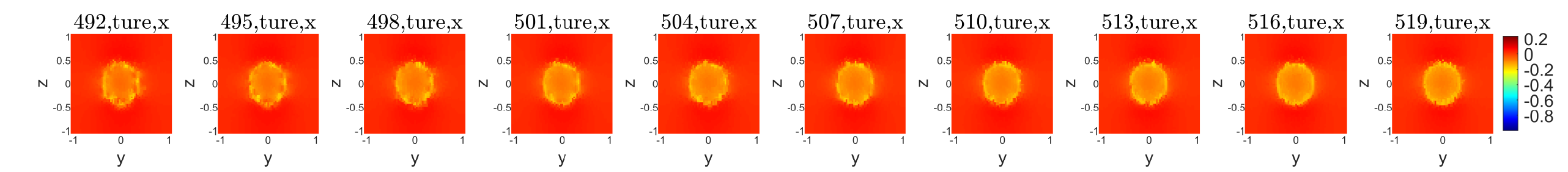}(a)\\
     \includegraphics[width=1\linewidth]{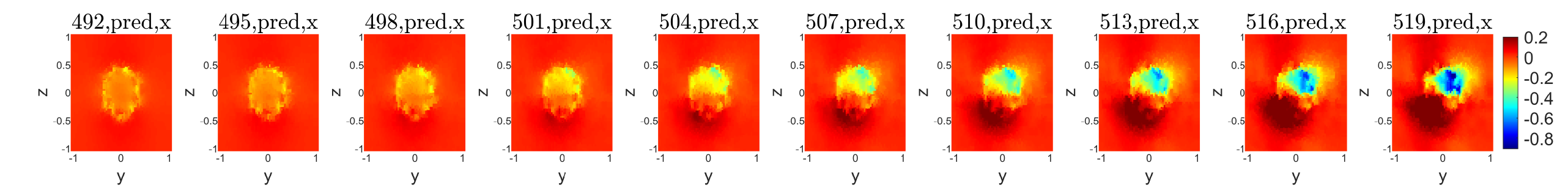}(b)\\
     \includegraphics[width=1\linewidth]{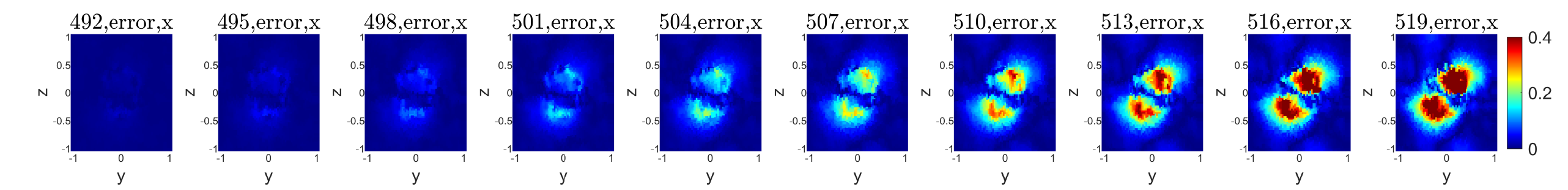}(c)\\
     \includegraphics[width=1\linewidth]{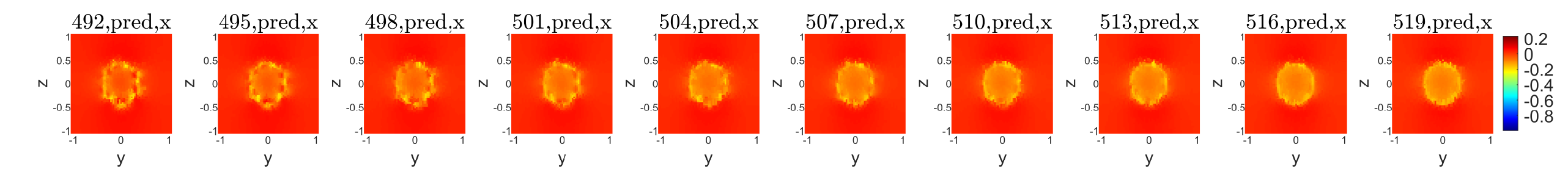}(d) \\
     \includegraphics[width=1\linewidth]{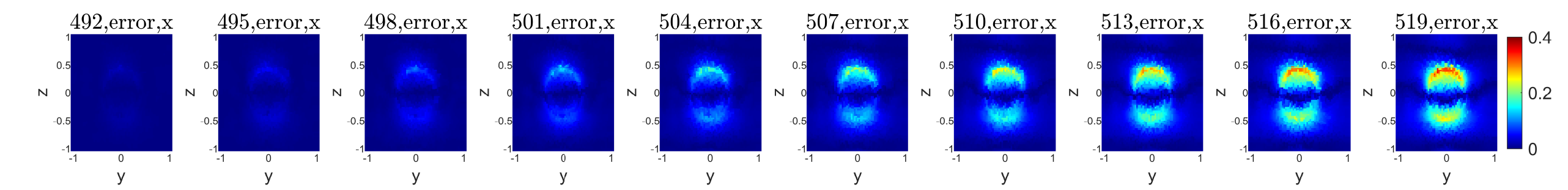}(e)\\
     \includegraphics[width=1\linewidth]{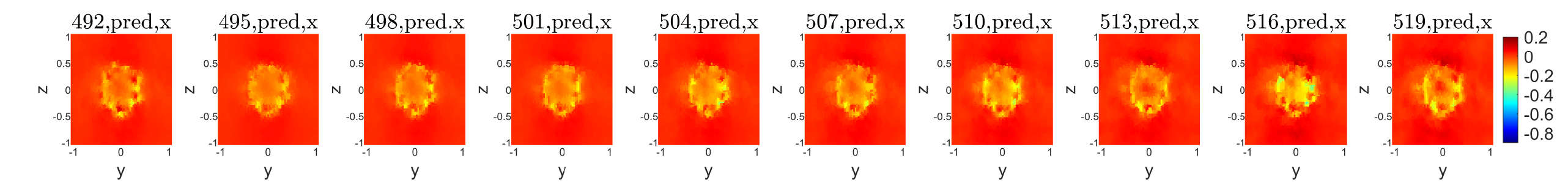}(f)\\
     \includegraphics[width=1\linewidth]{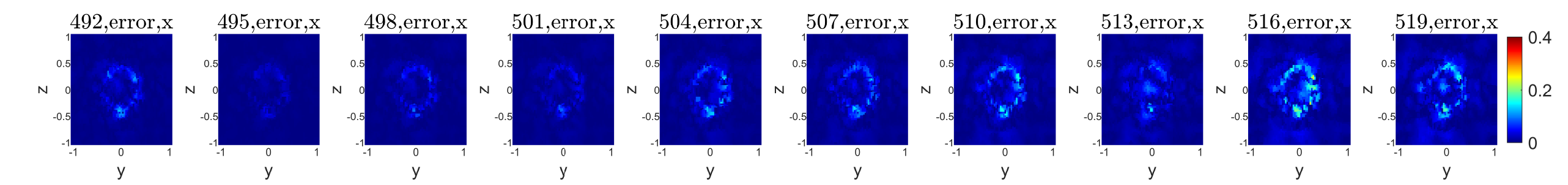}(g)\\
     \includegraphics[width=1\linewidth]{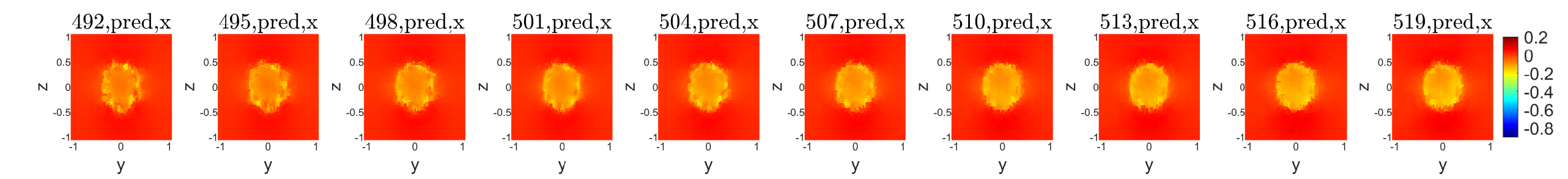}(h)\\
     \includegraphics[width=1\linewidth]{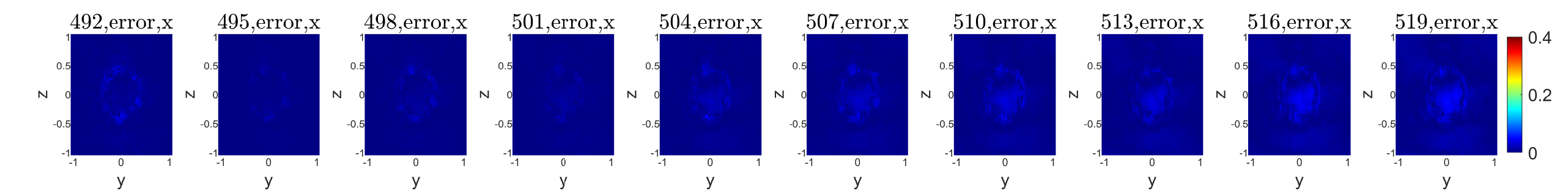}(l)\\    
     \caption{The predicted results and errors of the four FNOs with mix-set as input are shown. (a) represents the solution obtained by IBM. (b), (d), (f), and (h) represent the predicted solutions of FNO1, FNO2, FNO-3, and FNO-4, respectively. (c), (e), (g), and (l) represent the errors of FNO-1, FNO-2, FNO-3, and FNO-4, respectively.} 
     \label{exp482_491} 
 \end{figure}

\begin{figure}[h]
     \centering    
     \includegraphics[width=1\linewidth]{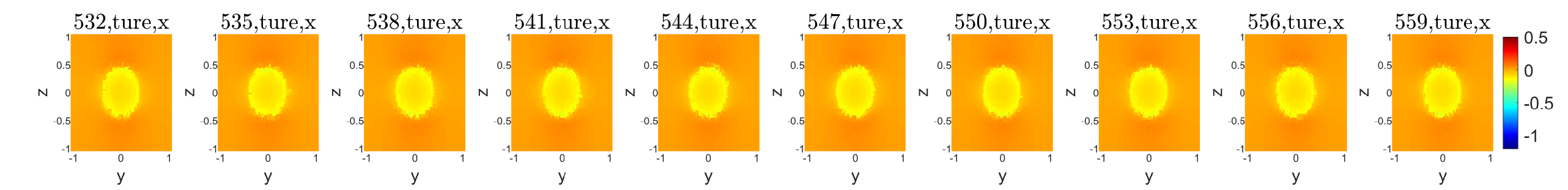}(a)\\
     \includegraphics[width=1\linewidth]{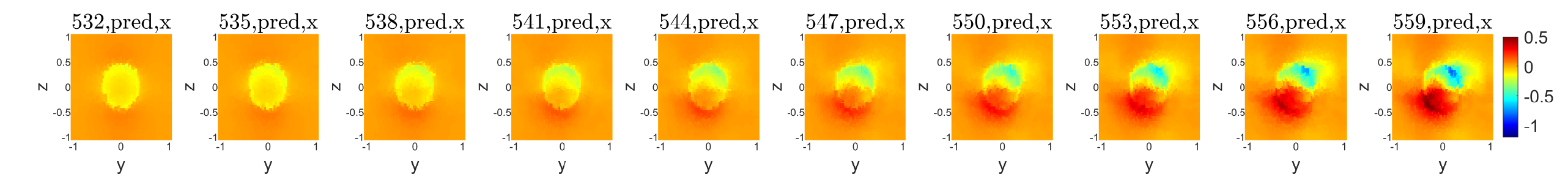}(b)\\
     \includegraphics[width=1\linewidth]{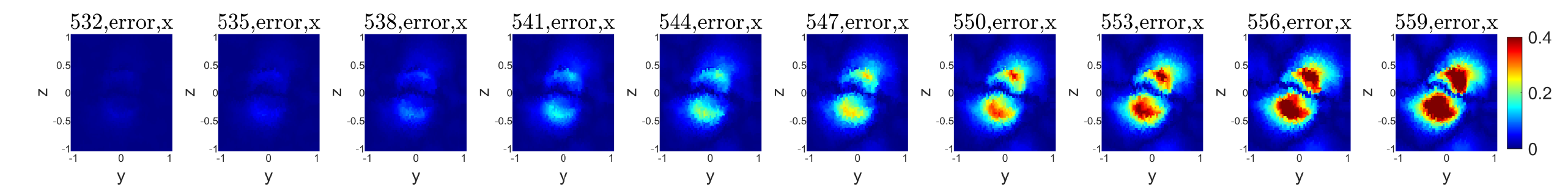}(c)\\
     \includegraphics[width=1\linewidth]{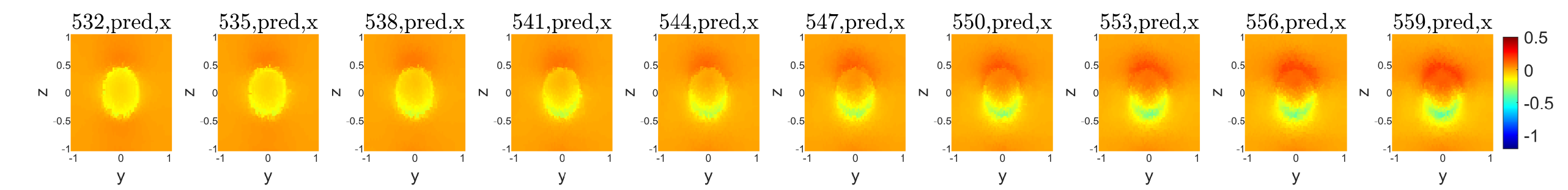}(d) \\
     \includegraphics[width=1\linewidth]{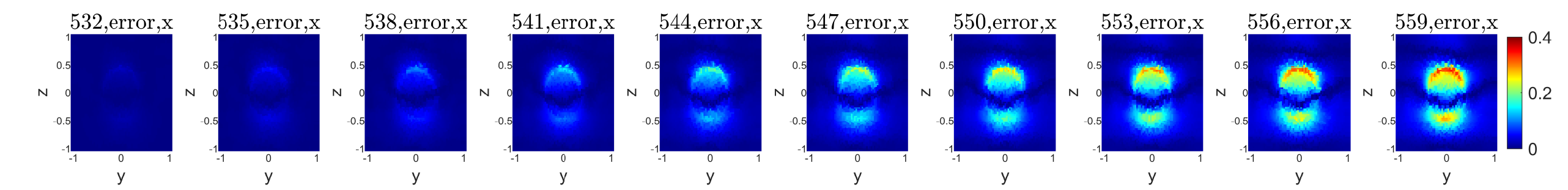}(e)\\
     \includegraphics[width=1\linewidth]{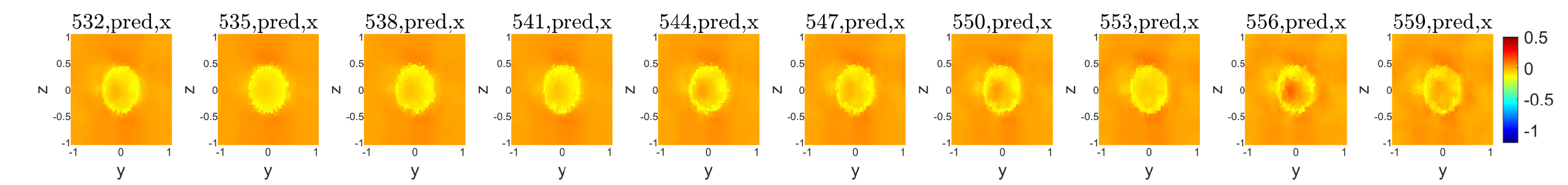}(f)\\
     \includegraphics[width=1\linewidth]{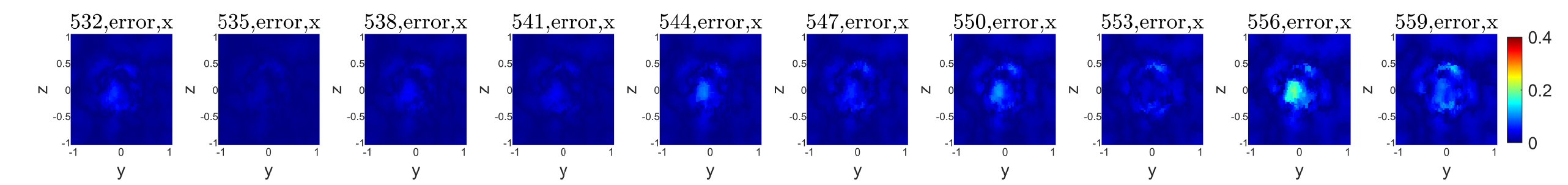}(g)\\
     \includegraphics[width=1\linewidth]{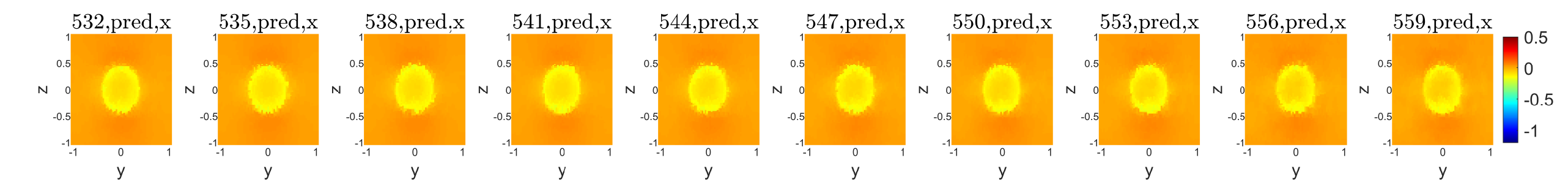}(h)\\
     \includegraphics[width=1\linewidth]{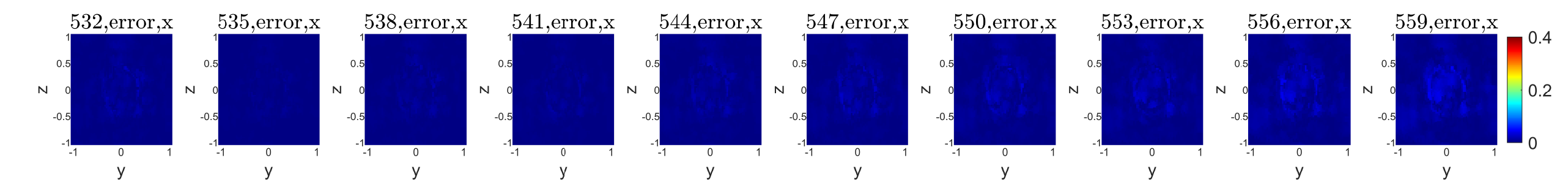}(l)\\    
     \caption{The predicted results and errors of the four FNOs with extra-set as input are shown. (a) represents the solution obtained by IBM. (b), (d), (f), and (h) represent the predicted solutions of FNO-1, FNO-2, FNO-3, and FNO-4, respectively. (c), (e), (g), and (l) represent the errors of FNO-1, FNO-2, FNO-3, and FNO-4, respectively.} 
     \label{exp522_531} 
 \end{figure}

The predicted results and errors of experiments, using inter-set, mix-set, and extra-set as input, are shown in Figs. \ref{exp442_451}, \ref{exp482_491}, and \ref{exp522_531}, respectively. In each subplot of the figures: (a) represents the IBM solution; (b), (d), (f), and (h) represent the predicted solutions of FNO-1, FNO-2, FNO-3, and FNO-4, respectively; (c), (e), (g), and (l) represent the errors of FNO-1, FNO-2, FNO-3, and FNO-4, respectively.

\begin{figure}[h]
     \centering
     \includegraphics[width=1\linewidth]{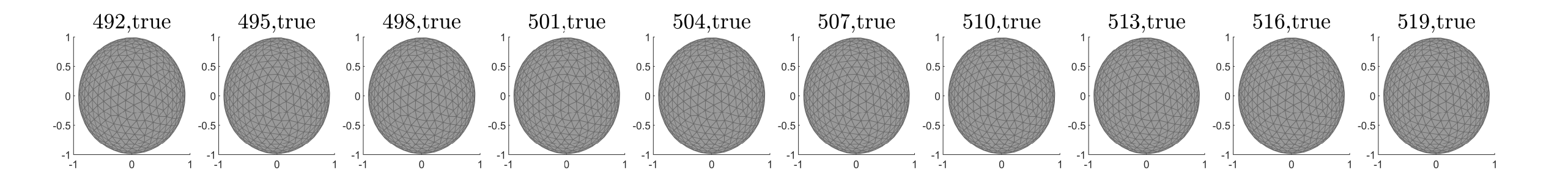}(a)\\
     \includegraphics[width=1\linewidth]{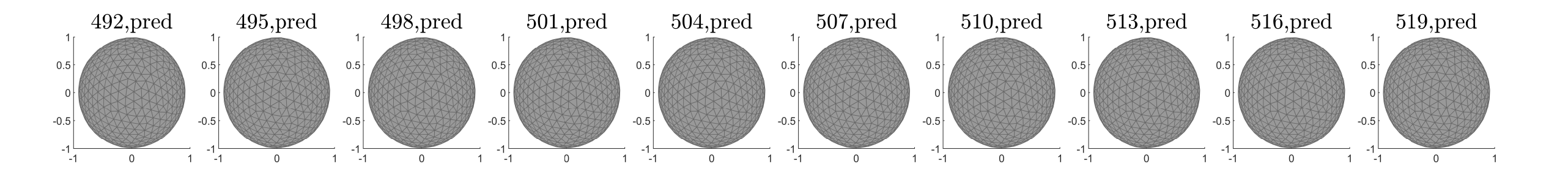}(b)\\
     \includegraphics[width=1\linewidth]{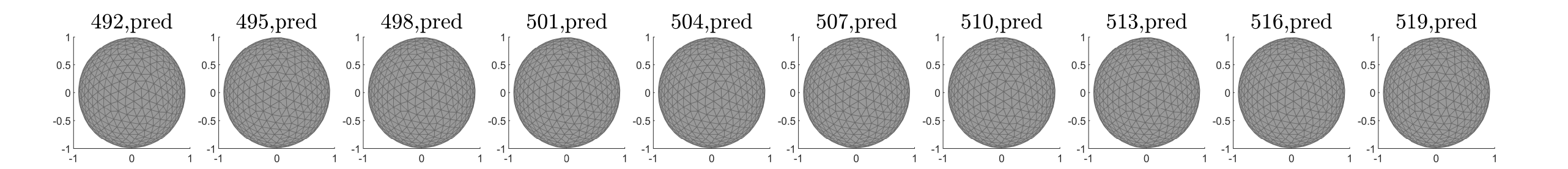}(c)
     \caption{The evolution of the vesicle at the last 30 moments with mix-set as input. (a) represents the solution obtained by IBM. (b) represents the solution obtained by using the fluid velocity predicted by FNO-3. (c) represents the solution obtained by using the predictions by FNO-4.}
     \label{v_492}
 \end{figure}

\begin{figure}[h]
     \centering
     \includegraphics[width=1\linewidth]{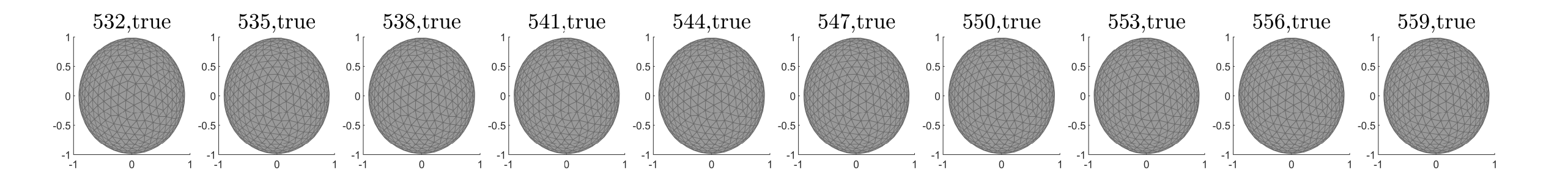}(a)\\
     \includegraphics[width=1\linewidth]{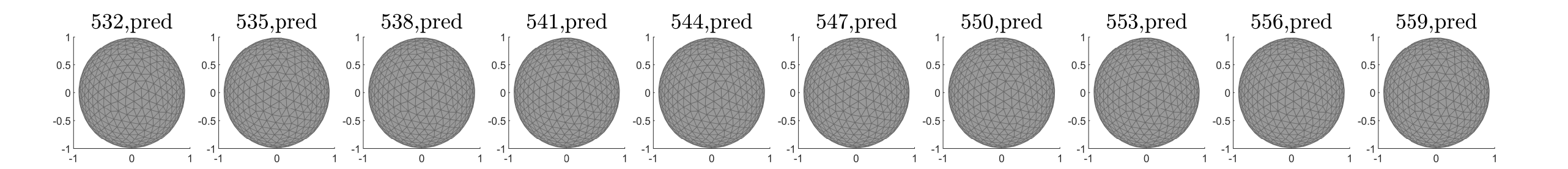}(b)\\
     \includegraphics[width=1\linewidth]{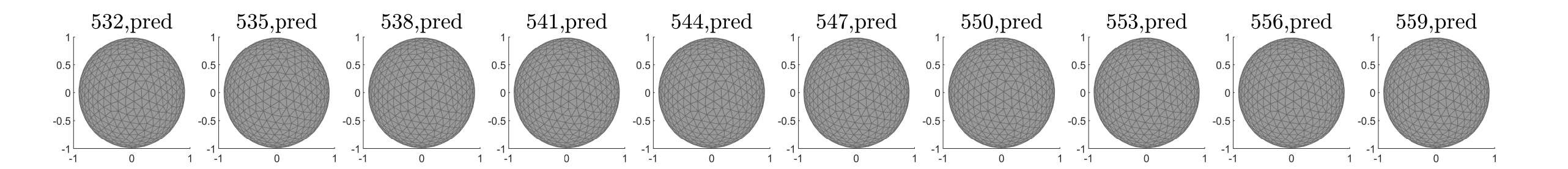}(c)
     \caption{The evolution of the vesicle at the last 30 moments with extra-set as input. (a) represents the solution obtained by IBM. (b) represents the solution obtained by using the fluid velocity predicted by FNO-3. (c) represents the solution obtained by using the predictions by FNO-4.}
     \label{v_532}
 \end{figure}

 The evolution of a vesicle at the last 30 moments with mix-set and extra-set as input are shown in Fig. \ref{v_492} and Fig. \ref{v_532}, respectively. Subfigure (a) represents the solution obtained by transferring the fluid velocity predicted by IBM. Subfigure (b) represents the solution obtained by transferring the prediction of FNO-3, and subfigure (c) represents the solution obtained by transferring the prediction of FNO-4.

\bibliography{cite}

\end{document}